\documentclass[11pt,twoside]{amsart}
\usepackage{amsfonts}
\usepackage{epsfig,graphics}
\usepackage{amssymb}
\usepackage{amscd}

\newtheorem{theoreme}{Theorem}[section]

\newtheorem{proposition}[theoreme]{Proposition}
       
\newtheorem{definition}[theoreme]{Definition}

\newtheorem{lemme}[theoreme]{Lemma}

\newtheorem{corollaire}[theoreme]{Corollary}

\newtheorem{question}[theoreme]{Question}
\newtheorem{remarque}[theoreme]{Remark}

\newcommand{\R}{{\bf R}}

\newcommand{\N}{{\bf N}}
\newcommand{\C}{{\bf C}}

\newcommand{\RP}{{\bf RP}}

\newcommand{\diag}{\text{diag }}
\newcommand{\add}{\text{ad }}

\newcommand{\Ein}{\text{Ein}}

\newcommand{\Conf}{\text{Conf }}

\newcommand{\Ad}{\text{Ad }}

\newcommand{\liei}{{\mathfrak{I}}}
\newcommand{\liez}{{\mathfrak{z}}}
\newcommand{\lien}{{\mathfrak{n}}}

\newcommand{\lieg}{{\mathfrak{g}}}
\newcommand{\lieh}{{\mathfrak{h}}}
\newcommand{\liep}{{\mathfrak{p}}}

\newcommand{\lier}{{\mathfrak{r}}}
\newcommand{\oo}{{\mathfrak{o}}}

\newenvironment{preuve}{\medskip \noindent {\bf Proof: }}
   {$\diamondsuit$ }

\begin{document}
\pagenumbering{arabic}
\title[Local dynamics]{Local dynamics  of conformal vector fields}
\author{Charles Frances}
\date{\today}
\address{Charles Frances. 
Laboratoire de
Math\'ematiques, 
Universit\'e Paris-Sud. 
91405 ORSAY Cedex.}
\email{Charles.Frances@math.u-psud.fr}
\maketitle

\begin{abstract}
We study pseudo-Riemannian conformal vector fields in the neighborhood of a singularity. For Riemannian manifolds, we prove that if a  conformal vector field vanishing at a point $x_0$ is not Killing   for a metric in the conformal class, then  a neighborhood of the singularity $x_0$ is conformally flat. \end{abstract}

\section{Introduction}
A smooth vector field on a pseudo-Riemannian manifold $(M,g)$ is {\it conformal}  whenever its local flow leaves invariant the conformal class $[g]:=\{ e^{\sigma}g \ | \ \sigma \in C^{\infty}(M) \}.$  Otherwise stated, $X$ is a conformal vector field when there exists $\psi \in C^{\infty}(M)$ such that  $ {\mathcal L}_Xg=\psi g.$ 

Locally, on a  pseudo-Riemannian manifold of type $(p,q)$ and dimension $n=p+q$, $n \geq 3$, the Lie algebra of conformal vector fields is finite dimensional, of dimension smaller or equal to $\frac{(n+2)(n+1)}{2}$. For every type $(p,q)$, there exists a compact, conformally flat manifold, for which this maximal dimension is attained, not only locally but also for the Lie algebra of global conformal vector fields. This manifold, called Einstein's universe $\Ein^{p,q}$, is obtained projectivizing in ${\bf RP}^{n+1}$ the null cone of a type-$(p+1,q+1)$  quadratic form  (see section \ref{sec.univers}).  Any $1$-parameter subgroup of ${\mbox O}(p+1,q+1)$ acting on $\Ein^{p,q}$  defines a conformal vector field. Such fields on $\Ein^{p,q}$ will be called {\it  M\"obius vector fields}. They are a handy family of model fields, that one can study explicitely. 

Studying conformal vector fields on a  general manifold $(M,g)$ is a more difficult problem. One would like to understand completely the behavior of such fields around a singularity, and if possible to exhibit a complete family of  {\it normal forms}. Though this problem has motivated quite a lot of works, mainly in Riemannian ($p=0$) or Lorentzian ($p=1$) geometry  (see, among other references  \cite{beig}, \cite{capocci1}, \cite{capocci2}, \cite{frances1}, \cite{kuehnel}, \cite{kuehnel2}, \cite{kuehnel4}, \cite{steller}...), the local understanding of conformal vector fields is still not fully satisfactory.   

A relevant notion, when studying conformal vector fields is that of  {\it essentiality}. If $X$ is a conformal vector field on $(M,g)$, vanishing at $x_0 \in M$, one says $X$ is {\it inessential} in the neighborhood of  $x_0$ whenever there is an open set  $U$ containing  $x_0$,  and a smooth function  $\sigma : U \to \R$ such that $X$ is a {\it Killing field} for  $e^{\sigma}g$ on $U$; this means ${\mathcal L}_X(e^{\sigma}g)=0$.  In this case, the exponential map of $e^{\sigma}g$ allows to linearize $X$ near $x_0$: indeed, it conjugates the action of the local flow $\{ \phi_X^t \}$ generated by $X$, and that of the flow of the differential  $\{D_{x_0}\phi_X^t\}$, on suitable neighborhoods of  $x_0$ and   $0 \in T_{x_0}M$ respectively.  Hence, inessential conformal vector fields are easy to understand.
Now, what can be said about {\it essential} vector fields, {\it i.e} those which do not preserve any metric in the conformal class $[g]$, on any neighborhood of the singularity $x_0$?  

Our aim in this paper is to use the canonical Cartan connection associated to a pseudo-Riemannian conformal structure to study conformal vector fields.  This is not really a new idea, since it  already appeared, for instance, in  \cite{aleks} or \cite{capocci1}, but we would like to sistemize it.

The upshot is that, as soon as we are dealing with manifolds of dimension at least three,  the canonical Cartan connection associated to a conformal structure allows to define  a distinguished class of curves: the conformal geodesics. This family is preserved by any local conformal transformation.    Next, a base point $o$ being fixed once for all on the model space $\Ein^{p,q}$, it is possible to associate naturally to $X$ and the singularity $x_0$ a M\"obius field $X_h$, vanishing at $o$: this is the {\it holonomy vector field} of $X$ at $x_0$.  This field  $X_h$ generates a $1$-parameter group $\{h^t\}: =\{\exp({tX_h}) \}\subset \mbox{O}(p+1,q+1)$, fixing $o$, and called the  {\it holonomy flow}. A fundamental fact is that the action of the local flow $\{ \phi_X^t\}$ on conformal geodesic segments emanating from $x_0$ is actually {\it conjugated} to the action of $\{h^t\}$ on the geodesic segments emanating from $o$.  Unfortunately, unlike what is happening in the isometric case, one can not infer directly a local conjugacy between   $X$ and $X_h$. This is basically because there are ``too many" conformal geodesics   for the ``conformal exponential map" being a local diffeomorphism in a neighborhood of the origin.  Yet, we will show that much of the local  properties of $X$ around $x_0$ are encoded in its holonomy field $X_h$, and our aim  is to develop this idea of a dictionnary between a conformal field and its holonomy, and to give evidences that there should be a positive answer to  the following:  

\begin{question}
\label{quest.conjugaison}
In dimension $n \geq 3$, is a pseudo-Riemannian conformal vector field $X$ always locally conjugated in a neighborhood of a singularity $x_0$ to its holonomy vector field  $X_h$ (in particular to a M\"obius vector field)?
\end{question}

The formulation is deliberately vague concerning the regularity of the conjugacy, if it exists. We may have only a $C^0$-conjugacy between orbits of the local flows.

One of the main results of the article, is to get a complete description of  Riemannian conformal vector fields around a singularity.  In particular, we answer affirmatively to question  \ref{quest.conjugaison} in this case (we will consider this question in Lorentzian signature in \cite{champs-karin}).  More precisely, we prove:  
  

\begin{theoreme}
\label{thm.ferrand-champs}
Let $(M,g)$ be a smooth Riemannian manifold of dimension  $n \geq 3$,  endowed with a smooth conformal vector field $X$.  We assume that $X$ vanishes at $x_0 \in M$. Then :

\begin{enumerate}
\item{Either there is a neighborhood  $U$ of $x_0$ on which  $X$  is complete,  and generates a flow which is relatively compact  in  $\Conf (U)$.  In this case  $X$ is linearizable around  $x_0$, and is inessential on $U$.}
\item{If we are not in the previous case, there is a conformally flat neighborhood  $U$ of $x_0$. The field $X$ is essential on each neighborhood of $x_0$.}
\end{enumerate}
In both cases, $X$ is smoothly conjugated in the neighborhood of $x_0$ to its holonomy field $X_h$.
\end{theoreme} 

This theorem is actually a particular case of theorem \ref{thm.ferrand-algebres}, which deals not only with single vector fields, but with any subalgebra of conformal vector fields.

If one adds a completeness assumption on the field $X$,  theorem \ref{thm.ferrand-champs} reduces to a particular case of the so-called Obata-Ferrand's theorem (see \cite{obata}, \cite{ferrand}, \cite{schoen}, \cite{frances-ferrand}). This theorem states that a group of {\it global} conformal transformations on a  Riemannian manifold   $(M,g)$ is always inessential  ({\it i.e} it acts isometrically for a metric in the conformal class  $[g]$) except maybe if $(M,[g])$ is conformally equivalent to the standard sphere or the Euclidean space.
 The (independant) proofs of  \cite{ferrand} and \cite{schoen} rely on the use of global conformal invariants, and thus do not seem to carry over into  the study of conformal vector fields. 
In \cite{aleks}, D. Alekseevskii proposed a proof of theorem \ref{thm.ferrand-champs} in the case of complete vector fields, but some of his arguments turned out to be flawed.  Nevertheless, Section 2 of \cite{aleks}  develops ideas which are very similar to our notion of holonomy  (see also \cite{ferrand-correction} for some complements about Aleksevskii's work).

\subsection{Organisation of the paper}
Sections \ref{sec.geometrie} and \ref{sec.geodesiques}  introduce some material about Einstein's universe, the interpretation of conformal structures of dimension  $\geq 3$ as Cartan geometries, and the definition of conformal geodesics.  

Then comes the heart of the paper, with the notion of holonomy of a vector field, introduced in section \ref{sec.holonomy}. We  give a characterization  of local inessentiality, as well as linearizability, in terms of holonomy  (propositions \ref{prop.lineaire} and \ref{prop.essentiel}).

The main difficulty in theorem  \ref{thm.ferrand-champs} is that the vector fields we are considering are not assumed to be complete: all the orbits of the local flow  $\{\phi_X^t\}$ (except that of  $x_0$) could blow-up in finite time.  It is thus crucial to exhibit conditions, for any type-$(p,q)$ manifold, ensuring that  a conformal vector field $X$ admits ``a lot" (actually a nonempty open set) of {\it semi-complete} orbits.  By semi-complete, we mean defined on  $]-\infty,0]$ or $[0,\infty[$. This will be done  in sections   \ref{sec.dynamique-holonomy} and  \ref{sec.semi-complet}.  

When one is able to show that some orbits are semi-complete, then one can study asymptotic properties of those orbits, and using dynamical arguments, get geometrical informations such as, for instance, conformal flatness.  This is the basic idea behind theorem  \ref{thm.ferrand-algebres}, which is proved  in the last section. 

\section{Conformal structures and  Cartan geometries}
\label{sec.geometrie}
 
\subsection{Einstein's universe}
\label{sec.univers}
In a lot of geometrical contexts, there exists among all possible structures,  a distinguished object characterized by its great amount of symmetries.  For type-$(p,q)$ conformal structures, this object is Einstein's universe  $\Ein^{p,q}$, of which we give now a brief description  (see \cite{primer} or \cite{charlesthese} for more details).  In the whole paper, we assume that $p \leq q$ (in case $q>p$, one should switch the words ``time" and ``space" in all the statements).

Let us call $\R^{p+1,q+1}$  the space  $\R^{p+q+2}$ endowed with the quadratic form:
$$ Q^{p+1,q+1}(x):=2x_0x_{p+q+1}+...+2x_px_{q+1}+\Sigma_{p+1}^q x_i^2$$
The null cone of $Q^{p+1,q+1}$ is denoted ${\mathcal N}^{p+1,q+1}$.  When restricted to ${\mathcal N}^{p+1,q+1}$, the form   $Q^{p+1,q+1}$ yields a degenerate metric, the ($1$-dimensional) kernel of which is  tangent to the cone. Thus, the projectivization  ${\bf P}({\mathcal N}^{p+1,q+1} \setminus \{ 0 \})$ is a smooth submanifold of $\RP^{p+q+1}$, naturally endowed with a conformal class of type-$(p,q)$ metrics.  
One calls {\it Einstein's universe} of type  $(p,q)$, denoted by  $\Ein^{p,q}$, this compact manifold ${\bf P}({\mathcal N}^{p+1,q+1} \setminus \{  0  \})$ endowed with the conformal structure described above.

Notice that the space  $\Ein^{0,q}$ is merely the sphere ${\bf S}^q$ endowed with the conformal class of the round metric $g_{{\bf S}^q}$.
For  $p \geq 1$,  the product ${\bf S}^p \times {\bf S}^q$, endowed with the conformal class of the product metric  $-g_{{ \bf S}^p} \oplus g_{{\bf S}^q}$ is a double cover of de $\Ein^{p,q}$.

Let $\mbox{O}(p+1,q+1)$ be the group of linear transformations preserving $Q^{p+1,q+1}$.  Clearly, the natural action of  $\mbox{PO}(p+1,q+1)$ on $\Ein^{p,q}$ preserves the conformal class of $\Ein^{p,q}$, and actually  $\mbox{PO}(p+1,q+1)$ turns out to be the whole group of conformal transformations of $\Ein^{p,q}$  (this is Liouville's theorem, see \cite{charlesthese} or \cite{schot}, and references in  \cite{kuehnel3}).

In the whole paper, we will call $o$ the point of $\Ein^{p,q}$ corresponding to  $[e_0]$.  Its stabilizer  $P \subset \mbox{PO}(p+1,q+1)$ is a parabolic subgroup isomorphic to the semi-direct product  $(\R_+^* \times \mbox{O}(p,q)) \ltimes \R^{p,q}=\Conf(\R^{p,q})$.  From the conformal point of view, $\Ein^{p,q}$  is the homogeneous space $\mbox{PO}(p+1,q+1)/P$.  

Let $j : \R^{p,q} \to \Ein^{p,q}$ be the map given in projective coordinates on ${\bf RP}^{p+1,q+1}$ by:
$$ j : \left( \begin{array}{c} x_1\\
\vdots\\
x_n\\
\end{array} \right) \mapsto \left[ \begin{array}{c}-\frac{Q^{p,q}(x)}{2}\\
x_1\\
 \vdots\\
  x_n \\
   1 \\
   \end{array}
   \right]$$
The map $j$ is a conformal embedding from  type-$(p,q)$ Minkowski's space  $\R^{p,q}$ onto a dense open subset of  $\Ein^{p,q}$.  This map is the {\it stereographic  projection} of pole $o$. The image $j(\R^{p,q})$ is the complementary in $\Ein^{p,q}$ of the  {\it lightcone } of vertex $o$, namely the set of all lightlike geodesics emanating from $o$.  By a slight abuse of  language, we will often identify in the sequel $\R^{p,q}$ and its image $j(\R^{p,q})$, namely we will see Minkowski's space as an open subset of  Einstein's  universe.  Notice that $j$ conjugates the action of  $P$ on $j(\R^{p,q})$ and the affine action of $P=(\R_+^* \times \mbox{O}(p,q)) \ltimes \R^{p,q}$ on $\R^{p,q}$.

The basepoint $o$ is not in $j(\R^{p,q})$; it is located ``at infinity". It is thus convenient to introduce a second conformal chart:
$$ j^o: \left( \begin{array}{c}
x_1\\
\vdots\\
x_n \\
\end{array}
\right) \mapsto \left[ \begin{array}{c}
1 \\
 x_1 \\
  \vdots \\
  x_n \\
   -\frac{Q^{p,q}(x)}{2}\\
   \end{array}
   \right]$$
The map  $j^o$ is a conformal diffeomorphism from  $\R^{p,q}$ to an open subset of  $\Ein^{p,q}$ containing $o$. We denote $\R_o^{p,q}$ this open set and we say that  $j^o$ is the  ``chart at infinity ".    

Finally, let us notice that  $\R_o^{p,q} \cap j(\R^{p,q})$ is simply the image by  $j$ of $\R^{p,q}$ with its null cone (with vertex $0$) removed.

\subsection{The  Lie algebra $\oo(p+1,q+1)$}
\label{sec.algebre}
The  Lie algebra $\oo(p+1,q+1)$ comprises all matrices $X$ of size 
$(p+q+2) \times (p+q+2)$ satisfying the indentity: 
$$ X^t J_{p+1,q+1} + J_{p+1,q+1} X = 0$$
Here,  $J_{p+1,q+1}$ is the matrix of the quadratic form  $Q^{p+1,q+1}$ expressed in the base $(e_0, \ldots ,e_{n+1})$.

The algebra $\oo(p+1,q+1)$writes as a sum $\lien^- \oplus \lier \oplus \lien^+$, where:

$$ \lier =  \left\{ \left( \begin{array}{ccc}
a &  & 0 \\
  & M   &  \\
  &    & -a
\end{array} \right) \ :    
\qquad 
\begin{array}{c}
 a \in \R  \\
  M \in \oo(p,q) \\
    \end{array} 
\right\}
$$

$$ \lien^+= \left\{ \left( \begin{array}{ccc}
0& -x^t.J_{p,q}   &  0\\
  & 0  & x \\
  &   & 0
\end{array} \right) \ :  
\qquad 
\begin{array}{c}
  x\in \R^{p,q} 
\end{array}
\right\}
$$

$$ \lien^- = \left\{ \left( \begin{array}{ccc}
0&    &  \\
 x & 0  &  \\
 0 & -x^t.J_{p,q}  & 0
\end{array} \right) \ :  
\qquad 
\begin{array}{c}
  x\in \R^{p,q} 
\end{array}
\right\}
$$

In $\oo(p+1,q+1)$, one calls ${\mathfrak a}$ the algebra comprising  the matrices:
$${\mathfrak a} = \left\{ \left( \begin{array}{ccccccc}
\alpha_1&    & & & & &  \\ $$
 & \ddots & & & & &  \\
 & & \alpha_{p+1} & & &  &\\
 & & & I_{q-p} & & &  \\
 & & & & -\alpha_{p+1}& & \\
 & & & & & \ddots & \\
 & & & & & & -\alpha_{1} \\
 \end{array}
\right)  :  
\qquad 
\begin{array}{c}
  \alpha_1, \ldots ,\alpha_{p+1} \in \R 
\end{array} \right\}.
$$

The closed subgroup of  $P$ with Lie algebra ${\mathfrak a}$ is denoted by  $A$. Let us call  ${\mathfrak a}^+$ the subset of  ${\mathfrak a}$ corresponding to  $\alpha_1, \ldots , \alpha_{p+1} \geq 0$, and  $A^+:= \exp({{\mathfrak a}^+}).$  

\subsection{The parabolic subgroup $P$}
\label{sec.matriciel}
The group  $P$ is the stabilizer in ${\mbox O}(p+1,q+1)$ of the point $o=[e_0]$.  As already seen, it is isomorphic to the  semi-direct product $(\R_+^* \times \mbox{O}(p,q)) \ltimes \R^{p,q}$.
Using the chart $j$, we will often see the elements of  $P$ as affine transformations $A+T$, where $A \in \R_+^* \times \mbox{O}(p,q)$ is the linear part, and  $T \in \R^{p,q}$  is the  translation factor.  

As an element of  ${\mbox O}(p+1,q+1)$, a   translation of vector  $v \in \R^{p,q}$ writes:
$$ n^+(v):=\left( \begin{array}{ccc}
1& -v^t.J_{p,q}   &  -\frac{Q^{p,q}(v)}{2}\\
  & 1  & v \\
  &   & 1
\end{array} \right).  
 $$
The set of translations constitutes a group $N^+$, with Lie algebra  $\lien^+$.  The map  $n^+ : \R^{p,q} \to N^+$ is a group isomorphism.

Seen in ${\mbox O}(p+1,q+1)$,  an element $\lambda A \in \R_+^* \times \mbox{O}(p,q)$ writes as:
$$ \left( \begin{array}{ccc}
\lambda & 0 & 0 \\
0  & A   & 0 \\
0  & 0   & \frac{1}{\lambda}
\end{array} \right)$$

\subsection{The problem of equivalence}
\label{sec.principe-equivalence}
Let  $G$ be a  Lie group, $P \subset G$ a closed subgroup, and  ${\bf X}=G/P$.  A {\it Cartan geometry} modelled on ${\bf X}$ is the data of a triple $(M,B,\omega)$, where:
\begin{enumerate}
\item{$M$ is a manifold having the same dimenson as  ${\bf X}$.}
\item{$\pi: B \to M$  is a $P$-principal bundle over $M$.}
\item{The form  $\omega$ is a $1$-form on $B$ taking values in   $\lieg$, and satisfying:
\begin{enumerate}
\item{ For every  $b \in B$, $\omega_b : T_bB \to \lieg$ is an isomorphism of vector spaces. }
\item{For every  $X \in \lieg$ and  $b \in B$, $\omega_b(\frac{d}{dt}_{|t =0}R_{\exp(tX)}.b)=X$.}
\item{For every  $p \in P$, $(R_p)^*\omega = (\Ad p^{-1}).\omega$.}
\end{enumerate}}
\end{enumerate}
Here, $R_p$ denotes the  right action by $p$ on $B$, and  $\exp$ is the  the exponential map of $G$. A $1$-form $\omega$ as above is called a  {\it Cartan connexion} on $B$.

One can see a Cartan geometry as a curved analogue of the {\it flat model } $({\bf X},G,\omega_G)$, where   $\omega_G$ is the Maurer-Cartan form on $G$.

Let us now take the example of the homogeneous space  ${\bf X}:=\Ein^{p,q}=\mbox{PO}(p+1,q+1)/P$, and of a  Cartan geometry $(M,B,\omega)$ modelled on  $\Ein^{p,q}$. Let us call $\pi:B \to M$ the bundle map. For every  $x \in M$ and $b \in B$ above  $x$, there exists a natural isomorphism:
$$ \iota_{b} : \lieg/\liep \to T_xM$$ defined by $\iota_{b}(\overline \xi):=D_b\pi(\omega_b^{-1}(\xi))$, where $\xi$ represents  the class  $\overline \xi \in \lieg/\liep$ in $\lieg$.  The isomorphism  $\iota_b$ satisfies the equivariance relation:
\begin{equation}
\label{equ.equivariance}
 \iota_{b.p^{-1}}((\Ad p).\overline \xi)= \iota_{b}(\overline \xi), \ \forall p \in P
\end{equation}

So, if  ${\mathcal C}$ denotes the unique conformal class of type-$(p,q)$ scalar products which are  $(\Ad P)$-invariant on $\lieg/\liep$, $\iota_b({\mathcal C})$ determines a type-$(p,q)$  conformal class on  $T_xM$, which does not depend on the choice of  $b \in B$ above  $x$.  In other words, a Cartan geometry $(M,B,\omega)$ modelled on $\Ein^{p,q}$ defines a conformal class  $[g]$ of type-$(p,q)$ metrics on $M$.   

If the bundle $B$ is fixed, there are {\it a priori} a lot of Cartan connections   $\omega$ defining the class  $[g]$ on  $M$.  The following theorem ensures that there is a suitable choice of normalizations which makes  $\omega$ unique (this is the conformal analogue of the  Levi-Civita connection).  Precisely, see \cite[ch 7]{sharpe}, \cite{kobayashi}, one can state:

\begin{theoreme}[E. Cartan]
\label{thm.principe-equivalence}
Let $(M,[g])$ be a type-$(p,q)$ pseudo-Riemannian conformal structure, $p+q \geq 3$. Then there exists a unique normal Cartan geometry $(M,B,\omega)$, modelled on  $\Ein^{p,q}$, defining the conformal structure  $(M,[g])$ by the process described above. In particular, any local conformal diffeomorphism $\phi$ on  $M$ lifts to a local automorphism of the bundle $B$ (also denoted $\phi$) preserving $\omega$.
\end{theoreme}

In what follows, we will call the triple  $(M,B,\omega)$ given by theorem  \ref{thm.principe-equivalence}, the  {\it normal  Cartan bundle} defined by   $(M,[g])$.



\section{Conformal geodesics}
\label{sec.geodesiques}
In all this section, $(M,[g])$ denotes a type-$(p,q)$ pseudo-Riemanniann conformal structure, of  dimension $n \geq 3$.  We see this conformal structure as a Cartan geometry modelled on  $\Ein^{p,q}$,  and we call $(M,B,\omega)$ the corresponding normal   Cartan bundle. In the following, we will often call $G$ the Lie group $\mbox{PO}(p+1,q+1)$. 
\subsection{Development of curves }
\label{sec.development-courbes}
Let  $b \in B$ and  $\hat \alpha : I \to B$ a $C^1$-curve such that $\hat \alpha(t_0)=b$.  One defines the   {\it development  of $\hat \alpha$} at  $b$, denoted ${\mathcal D}_b(\hat \alpha)$, as the unique curve   $\hat \beta : I \to G$  satisfying   $\hat \beta (t_0)=1_G$ and:
$$ \omega(\alpha^{\prime})=\omega_G(\beta^{\prime})$$

Let now be  $x \in M$, $b \in B$  in the fiber of  $x$, and $\alpha : I \to M$  a  $C^1$-curve such that   $\alpha(t_0)=x$,  then one defines its development at $x$ relatively to  $b$ as follows:
$${\mathcal D}_x^b(\alpha)(t):=\pi_G({\mathcal D}_b(\hat \alpha)(t)),$$
where  $\pi_G$ stands for the  projection from  $G$ onto  $\Ein^{p,q}=G/P$, and  $\hat \alpha$ is a lift of  $\alpha$ in  $B$, such that $\hat \alpha(t_0)=b$.
The definition of  ${\mathcal D}_x^b(\alpha)$ does not depend on the lift $\hat \alpha$ because if  $\lambda : I \to B$ and  $p : I \to P$ are two  $C^1$-curves, and if we call  $\gamma(t)=\lambda(t).p(t)$, then:
\begin{equation}
\label{equ.formule-Cartan}
\omega(\gamma^{\prime}(t))=(\Ad p(t))^{-1}.\omega(\lambda^{\prime}(t))+\omega_G(p^{\prime}(t))
\end{equation}
Otherwise stated, ${\mathcal D}_b(\gamma)(t)={\mathcal D}_b(\lambda)(t).p(t)$ and thus (see  {\cite[p. 208]{sharpe}}):
 $$\pi_G({\mathcal D}_b(\gamma)(t))=\pi_G({\mathcal D}_b(\lambda)(t).p(t))$$  

\subsection{Conformal exponential map}
  The data of  $Z$ in $\lieg=\oo(p+1,q+1)$ defines naturally a vector field  $\hat Z$ on $B$ by the  relation $\omega(\hat Z)=Z$.  If $Z \in \lieg$, we call $\psi_Z^t$ the local flow  generated on $B$ by the field $\hat Z$.  At each $b \in B$, we define ${\mathcal W}_b \subset \lieg$ as the set of $Z$ such that  $\psi_Z^t$ is defined for   $t \in [0,1]$ at $b$.  Then one defines the exponential map at  $b$:
  $$\exp(b, \ ) : {\mathcal W}_b \to B$$
   as:
  $$ \exp(b,Z):=\psi_{\hat Z}^1.b$$
  
  It is a standard fact that  ${\mathcal W}_b$ is a neighborhood of  $0$, and the map $\xi \mapsto \exp(b, \xi)$ determines a diffeomorphism from an open set ${\mathcal V}_b \subset {\mathcal W}_b$ containing $0$ onto a neighborhood of  $b$ in  $B$.
 
 One checks easily that if $x \in M$, $b \in B$ above $M$, $\xi \in \lieg$ and  $\alpha(s):=\pi(\exp(b,s\xi))$, then the development $\beta(s):={\mathcal D}_{x}^{b}(\alpha)(s)$ is given by  $\beta(s)=\pi_G(\exp({s \xi}))$.
  
 Let $f$ be a conformal transformation of  $M$. Then $f_*(\hat Z)=\hat Z$, and if  $p \in P$, $(R_p)_*(\hat Z)=\hat{Z_p}$ , where $Z_p:={ (\Ad p^{-1}).Z}$.  One infers the following  important equivariance property: 
 \begin{equation}
 \label{equ.equivariance-exponentielle}
  f(\exp(b,\xi)).p^{-1}=\exp(f(b).p^{-1}, (\Ad p).\xi)
  \end{equation}

\subsection{Conformal geodesics}  
\label{sec.geodesiques-conformes}
The Cartan connection and  the exponential map allow to define a distinguished class of curves on  $(M,g)$, namely the {\it conformal geodesics}.  We won't try here to make the link between the definition given below, and previous ones as, for instance, that given in  \cite{ferrand2} (see also \cite{cap}).  
We will call {\it  parametrized conformal geodesics} of $M$ through  $x \in M$, any curve $s \mapsto \pi(\exp(b,s \xi))$, where  $b \in B$ is in the fiber above  $x$, $\xi \in( \Ad P). \lien^-$, and  $s$ takes values in an open interval  $I$ containing  $0$. One says the geodesic is {\it timelike}  (resp. {\it spacelike}, resp. {\it lightlike}) if $\iota_b(\overline{\xi})$ is timelike (resp. spacelike, resp. lightlike), where $\overline{\xi}$ is the projection of  $\xi$ on $\lieg/\liep$.  Because of the relation $\exp(b,s\xi)=\exp(\exp(b,s_0),(s-s_0)\xi)$, the tangent vector to a timelike (resp. spacelike, resp. lightlike) conformal geodesic is everywhere timelike (resp. spacelike, resp. lightlike).
If $\alpha$  is a conformal geodesic defined on  $I$, and $s_0 \in I$, we will call  $[\alpha]$ the set  $\alpha([0,s])$, and say that  $[\alpha]$ is  a {\it conformal geodesic segment} (or shortly conformal segment) emanating from $x$.  

Let us describe a little bit more precisely the conformal geodesics in the model $\Ein^{p,q}$.
Following the formula given in  section \ref{sec.algebre} for elements of  $\lien^-$, a conformal geodesic writes in  projective coordinates:
$$ s \mapsto p.\left[ \begin{array}{c}1\\
 s w_1 \\
  \vdots \\
   sw_n \\
    -s^2\frac{Q^{p,q}(w)}{2}\\
    \end{array} \right]$$
where $w:=\left( \begin{array}{c} w_1\\
 \vdots \\
  w_n\\
  \end{array}
  \right) \in \R^{p,q}$ and $p \in P$.

The geodesic is timelike (resp.  spacelike)  if and only if $w$ is timelike (resp.  spacelike).  In this case, it writes $s \mapsto j(p.(-\frac{1}{s}\frac{2w}{Q^{p,q}(w)}))$.  Thus, any timelike (resp. spacelike) geodesic is of the form:
\begin{equation}
\label{eq.geod-espace}
s \mapsto j(\frac{1}{s}v+v_0)
\end{equation}
 where  $v \in \R^{p,q}$ is timelike (resp. spacelike), and $v_0 \in \R^{p,q}$.
 We see in particular that in the chart $j$, timelike or spacelike geodesic segments emanating from $o$ are half lines.
 
 In the chart at infinity $j^o$, a timelike (resp. spacelike) geodesic writes:
 \begin{equation}
\label{eq.geod-espace-infini}
 s \mapsto \frac{2s(v+sv_0)}{Q^{p,q}(v+sv_0)}
 \end{equation}
 
 Let us focus now on lightlike geodesics.  They write:
 $$ s \mapsto p.\left[
 \begin{array}{c}
 1\\
  s w_1 \\
   \vdots \\
    sw_n \\ 
    0\\
    \end{array} \right]$$
 From the matrices given in  section \ref{sec.matriciel} for elements  $p \in P$, it is easy to check that the lightlike geodesic segments emanating from $o$ and contained in  $R_o^{p,q}$ read in the chart $j^0$ as lightlike line segments emanating from $0$.  Their parametrization is of the form:
 $$ s \mapsto \left[
 \begin{array}{c}
 1\\
  \tau(s) w_1 \\
   \vdots \\
    \tau(s)w_n \\ 
    0\\
    \end{array} \right]$$
where $s \mapsto \tau(s)$ is an homographic transformation.

\subsection{Some auxiliary metrics}
\label{sec.metrique-auxiliaire}
To compare a curve and its development, we will introduce Riemannian metrics on $G$ and  $B$ as follows.
On $\R^{p,q}$, we denote  $<x,x>=x_1^2+ \ldots +x_n^2$ the  standard Euclidean product, and we call $||.||$ the  norm it defines. We pull back this scalar product on $\R_o^{p,q}$ by the map $(j^o)^{-1}$, what endows  $\R_o^{p,q}$ with a flat Riemannian metric  $\rho^o$.  
If $v:=(v_1, \ldots , v_n)^t$ (the transpose of $(v_1, \ldots , v_n)$), we set:
$$n^-(v):= \left( \begin{array}{ccc}
1&    &  \\
 -J_{p,q}.v^t & 1  &  \\
 -\frac{Q^{p,q}}{2} & v  & 1
\end{array} \right)$$
The morhism  $n^- : \R^n \to N^-$ is an  isomorphism and the relation $j^o(v+w)=n^-(v).j^o(w)$ shows that $N^-$ acts simply transitively on  $\R_o^{p,q}$,  and the metric $\rho^o$ is $N^-$-invariant.  Hence, it induces a left-invariant metric $\rho^-$ on  $N^-$.  We denote  $< \ , \ >_{\lien^-}$ the scalar product induced by this metric on  $\lien^-$, and  $||.||_{\lien^-}$ the associated norm.

In all the paper, if $r>0$, we will call  ${\mathcal B}(0,r)$ (resp. ${\mathcal S}(0,r)$) the open ball (resp. the sphere) centered at  $0$ and of radius $r$ in $\lien^-$ for the norm $||.||_{\lien^-}$.
The ball $B(o,r)$ of center $o$ and radius $r$ for the metric  $\rho^o$ is just  $\pi_G(\exp({{\mathcal B}(0,r)}))$.

Let us choose a scalar product  $< \ ,\ >_{\lieg}$ on  $\lieg$,  inducing $< \ , \ >_{\lien^-}$ on $\lien^-$.  The product  $< \ , \ >_{\lieg}$ allows to define a left-invariant Riemannian metric $\rho^G$ on  $G$.  Similarly, we define a Riemannian metric  $\rho^B$ on  $B$ by the formula:
$$ \rho_b^B(u,v):=<\omega_b(u),\omega_b(v)>_{\lieg}, \ b \in B, \ \ u,v \in T_bB$$

In what follows, if $\beta$ is a $C^1$-curve in  $\R_o^{p,q}$, we will denote $L^o(\beta)$  its length with respect to the metric $\rho^o$.

One checks immediately that if  $\alpha$  is a $C^1$-curve through $b \in B$, and if  $\beta:={\mathcal D}_b(\alpha)$, then the length of  $\alpha$ with respect to  $\rho^B$ is the length of $\beta$ with respect to  $\rho^G$.
Also, if  $\beta$ is a  $C^1$-curve of  $N^-$ through $1_G$, then its length relatively to  $\rho^G$ is $L^o(\pi_G(\beta))$, since by construction, the restriction of $\rho^G$ to $N^-$ is $\rho^-$.

\subsection{Degeneration properties}

A key idea of the paper will be to recover the dynamics of a sequence of local conformal transformations thanks to their action on the geodesic segments. For this, we will have to understand how a sequence of geodesic segments can degenerate. This will be done understanding the link between a sequence of geodesic segments, and the sequence of its developments. One now proves some lemmas in this direction, mostely inspired by   \cite{frances-ferrand}.   

 We first show that the curves through $x$ in  $M$, the development of which  are ``short" curves  in $\R_o^{p,q}$, are themselves short. 

\begin{lemme}
\label{lem.lemme-court}
 Let $x \in M$,  and $b \in B$ in the fiber of  $x$. For any neighborhood  $U$ of $x$,  there is a real $r_U$ such that if   $\alpha : \ ]a, b[ \to M$  is a $C^1$-curve through   $x$, and if  ${\mathcal D}_x^b(\alpha)$ is included in $\R_o^{p,q}$ with $L^o({\mathcal D}_x^b(\alpha)) \leq r_U$, then the curve  $\alpha$ is included in  $U$.
\end{lemme}

\begin{preuve}
let us denote by  $B(b,R)$ the ball of center $b$ and radius  $R$ for the metric  $\rho^B$. We first choose $r_U$ such that  $\pi(B(b,r_U)) \subset U$. Since ${\mathcal D}_x^b(\alpha)$ is included in  $\R_o^{p,q}$, there exists  $p : \ ]a,b[ \to P$, such that  $\beta(t):={\mathcal D}_b(\alpha)(t).p(t)$ is a curve of $N^-$through  $1_G$. The length of $\gamma(t):=\alpha(t).p(t)$ with respect to  $\rho^B$ is then equal to the length  $l$ of  $\beta$ with respect to  $\rho^G$. But  $l$ is nothing else than  $L^o({\mathcal D}_x^b(\alpha))$, as mentioned in the previous paragraph.  
 We conclude that $\gamma$ is included in $B(b,r_U)$,  and since $\gamma$ projects on  $\alpha$, we get that  $\alpha$ is included in  $U$.
\end{preuve}

We now prove:

\begin{proposition}
\label{prop.segment-court}
Let  $R>0$ and $[\alpha]$  a geodesic segment emanating from $o$ and contained in  $B(o,R)$, then:
$$ L^o([\alpha]) \leq 8n R$$
\end{proposition}

\begin{preuve}
we work in the chart at infinity  $j^o$, and we consider  a geodesic segment $[\alpha]$ emanating from  $o$,  and included in  $B(0,R)$ (the Euclidean ball of center  $0$ and radius  $R$).  If it is a lightlike geodesic segment, then as it was observed in section \ref{sec.geodesiques-conformes}, it is a line segment emanating from $0$ and included in  $B(0,R)$.   Its length $L^o$ is then at most $R$ and the proposition holds clearly in this case.

Now, if  $[\alpha]$ is timelike or spacelike, we know from  formula  (\ref{eq.geod-espace-infini}) that:
$$ \alpha(s)=\frac{2s(v+sv_0)}{Q^{p,q}(v+sv_0)}, \ s \in [0,s_0]$$
Thus, $\alpha(s):=(\alpha_1(s), \ldots , \alpha_n(s))=(\frac{P_1(s)}{Q_1(s)}, \ldots , \frac{P_n(s)}{Q_n(s)})$, where $P_1, \ldots , P_n,\\ Q_1, \ldots, Q_n $ are polynomials of degree at most  $2$. We infer that for  $i=1, \ldots, n$, the derivative  $\alpha_i^{\prime}(s)$ vanishes at most three times on  $[0,s_0]$.  Thus, the curve  $\alpha_i : [0,s_0] \to [-R, R]$ has length at most  $8R$.  Indeed, if  $\beta : [0,s_0] \to [-R, R]$ is a  $C^1$-curve the derivative of which vanishes at most  $m$ times, the length of this curve is at most  $2(m+1)R$.
We get finally the  inequality:
$$ \int_{0}^{s_0} \sqrt{(\alpha_1^{\prime}(s))^2+ \ldots + (\alpha_n^{\prime}(s))^2}ds \leq \int_{0}^{s_0}( |\alpha_1^{\prime}(s)|+ \ldots + |\alpha_n^{\prime}(s)|)ds \leq 8nR$$

This yields the  proposition.
\end{preuve}

\begin{corollaire}
\label{coro.tend-vers-0}
 Let  $x_0 \in M$, $b_0 \in B$ in the fiber of $x_0$.  If  $[\alpha_k]$  is a sequence of conformal geodesic segments emanating from  $x_0$, and if  ${\mathcal D}_{x_0}^{b_0}([\alpha_k])$ tends to  $o$, then $[\alpha_k] \to x_0$.
\end{corollaire}

\begin{preuve}
for every neighborhood  $U$ of $x_0$, the previous proposition yields an integer $K$ such that if $k \geq K$, then $L^o([\alpha_k]) < r_U $, where $r_U$ is given by lemma \ref{lem.lemme-court}.  This same lemma then says that  $[\alpha_k] \subset U$ for $k \geq K$.  
\end{preuve}

\section{Holonomy, linearizability and essentiality}
\label{sec.holonomy}
In all this section, we consider $(M,g)$ a type-$(p,q)$, smooth pseudo-Rieman\-nian manifold , $p+q \geq 3$.  The associated normal  Cartan  bundle is denoted by  $(M,B,\omega)$.  We assume there exists a smooth conformal vector field  $X$ on  $M$,  having a singularity $x_0 \in M$.  We are going to explain how to associate naturally to   $X$ a  M\"obius field $X_h$ : {\it its holonomy at $x_0$.}  Then, we will begin to establish a dictionnary between the properties of   $X$ in a neighborhood of  $x_0$  and those of  $X_h$ in a neighborhood of $o$, begining with two issues: linearizability and essentiality.

 \subsection{Holonomy morphism, holonomy vector field, and holonomy algebra}
\label{holonomy-champ}
Let $\liei_{x_0}$ be the Lie algebra of conformal vector fields of $M$ vanishing at $x_0$. 
 Let $X$ be a  vector field of $\liei_{x_0}$. 
 The local flow $\{ \phi_X^t \}$  lifts to $B$, yielding a vector field on $B$, still denoted $X$, satisfying ${\mathcal L}_X \omega =0$  
 (here, ${\mathcal L}_X$ is the Lie derivative with respect to $X$).  
 Thus, for each $b \in B$, we get a linear monomorphism $ s_b: \liei_{x_0} \to \lieg$ defined by:
 $$ s_b(X)=\omega_b(X(b)).$$ 
Now, if $X,Y \in \liei_{x_0}$, we have the following formula \cite[Lemma 2.1]{bfm}:
$$ \kappa_b(\overline{s_b(X)},\overline{s_b(Y)})=s_b([X,Y])+[s_b(X),s_b(Y)],$$
where $\kappa: B \to \text{Hom}(\Lambda^2(\lieg/\liep),\lieg)$ is the curvature function of $(M,B,\omega)$ (see \cite[Definition 3.22]{sharpe}), and $\overline{s_b(X)},\overline{s_b(Y)}$ are the projections on $\lieg/\liep$ of $s_b(X)$ and $s_b(Y)$.  We thus see that as soon as $b_0$ is a point in the fiber of $x_0$, $-s_{b_0}: \liei_{x_0} \to \liep$  is an embedding of Lie algebras, called {\it the holonomy morphism} at $b_0$. The image of $-s_{b_0}$ is into $\liep$, because $X$ vanishes at $x_0$ and its lift is then tangent to  the fiber of $x_0$.
The element $s_{b_0}(X) \in \liep$ defines a right-invariant vector field $X_h$ on $G=\mbox{O}(p+1,q+1)$: it will be called  the  {\it holonomy field} of  $X$ at $x_0$ (relatively to $b_0$).  The holonomy vector field is a conformal vector field of $\Ein^{p,q}$, hence a {\it M\"obius field}, vanishing at $o$.
The $1$-parameter group $\{h^t \}$ of $P$ defined by $h^t:=\exp(ts_{b_0}(X))$ will be called the {\it holonomy flow} of  $X$ at $x_0$ (relatively to $b_0$).
Let $\{ \phi_X^t \}$  be the local flow defined by $X$, lifted to $B$.  Then for every $t \in \R$, we have:
 $$\phi_X^t.b_0.h^{-t}=b_0$$
Using equation (\ref{equ.formule-Cartan}), it is easy to check  that if $\alpha : I \to M$ is a  $C^1$-curve through $x_0$,  and if the local flow  $\{\phi_X^t\}$ is defined at each point of  $\alpha$ for  $t \in [0,\delta]$, then the following  equivariance relation holds:
\begin{equation}
\label{equ.equivariance-holonomy}
 {\mathcal D}_{x_0}^{b_0}(\phi_X^t.\alpha)(s)=h^t.{\mathcal D}_{x_0}^{b_0}(\alpha)(s), \text{  for every } s \in I \text{ and } t \in [0,\delta].
\end{equation}

Also,  read into the chart $j$, the flow $\{h^t\}$ is affine on $\R^{p,q}$, thus can be written  $A^t+T_t$, where $A^t \in \R_+^* \times \mbox{O}(p,q)$ is the linear part, and  $T_t \in \R^{p,q}$ the translation part.  The relation:
$$ D_{x_0}\phi_X^t(\iota_b(\xi))=\iota_b((\Ad h^t).\xi)$$
identifies the linear part  $A^t$ with the differential $D_{x_0}\phi_X^t$ read in an orthonormal frame of $T_{x_0}M$.

Finally, if we are considering $\lieh \subset \liei_{x_0}$ a subalgebra, it is also possible to define its {\it holonomy algebra $\lieh_h:=s_{b_0}(\lieh)$}. Its {\it holonomy group} $H_h \subset P$ is the connected subgroup of $P$ having $\lieh_h$ as Lie algebra. 

The definition of the holonomy group (and of the holonomy algebra) of $\lieh$ at $x_0$ depends on the choice of  $b_0$ in the fiber of $x_0$. If  $b_0$ is replaced by  $b_0.p$, with $p \in P$, then $H_h$ (resp.  $\lieh_h$) is changed into $p.H_h.p^{-1}$ (resp. into $(\Ad p).\lieh_h$).  The holonomy of  $\lieh$ at $x_0$ is thus well defined up to conjugacy in $P$.  By a slight abuse of language, when we will speak about {\it the holonomy} of $\lieh$ at $x_0$, we will mean {\it  a representative of the conjugacy class} of all possible holonomies at  $x_0$.

\begin{remarque}
\label{rem.cas-plat}
In case  $\lieh$ is a Lie algebra of  conformal vector fields on a neighborhood  $V$ of  $o$ in $\Ein^{p,q}$, then  we merely have  $\lieh=\lieh_h$. This follows from the fact that the normal Cartan  bundle is the inverse image of $V$ in $(\Ein^{p,q},\mbox{O}(p+1,q+1),\omega_G)$,  and from Liouville's theorem.
\end{remarque}


\subsection{Holonomy and linearizability}
\label{sec.holonomy-linearisabilite}
%
The property for  $X$ to be linearizable in a neighborhood of  $x_0$ can be seen very easily on its holonomy flow. This is the content of the following proposition.  
\begin{proposition}
\label{prop.lineaire}
The field  $X$ is linearizable in a neighborhood of   $x_0$ if and only if its holonomy flow $\{h^t\} \subset \Conf(\R^{p,q})$ has a fixed point on  $\R^{p,q}$.  In this case, there is a $C^{\infty}$-diffeomorphism from a neighborhood of $x_0$ onto a neighborhood of   $o$ which conjugates $X$ and  $X_h$.
\end{proposition}

\begin{preuve}
we begin with the easiest part of the proposition. Let us assume that $\{h^t\}$ fixes a  point of  $\R^{p,q}$.  Considering  $b_0.p$ instead of  $b_0$ for a suitable  $p \in P$ ({\it i.e} conjugating  $\{h^t\}$ by $p$), we can assume that this fixed point is $0$, namely $\{ h^t \} \subset \R_+^* \times \mbox{O}(p,q)$. 
 We then pick $r>0$ small enough so that  $\xi \mapsto \pi(\exp(b_0,\xi))$ is a diffeomorphism, denoted $\psi$, from ${\mathcal B}(0,r)$ onto a neighborhood $U$ of $x_0$, and  $\xi \mapsto \pi_G(\exp({\xi}))$ is a diffeomorphism, denoted $\varphi$, from ${\mathcal B}(0,r)$ onto a neighborhood  $V$ of $o$. For $t$ near  $0$, and  $\xi \in {\mathcal B}(0,r)$, the  relation:
 $$ \phi_X^t.\exp(b_0,\xi).h^{-t}=\exp(b_0,(\Ad h^t).\xi)$$
just writes $\phi_X^t \circ \psi = \psi  \circ (\Ad h^t)$, and is  available on ${\mathcal B}(0,r)$.  Thus, $\psi  $ conjugates the local flow of  $X$ on $U$ to the linear flow  $(\Ad h^t)$ on ${\mathcal B}(0,r)$. Moreover, $\psi \circ \varphi^{-1}$ is a smooth diffeomorphism from $U$ onto $V$ which conjugates $X$ and $X_h$.

We are now going to prove the converse statement, namely that if $\{h^t\}$ does not fix a point in $\R^{p,q}$, then $X$ is not linearizable around $x_0$.
In the chart $j$, $\{h^t\}$ is a flow of affine conformal transformations of $\R^{p,q}$, that we denote $h^t=e^{\lambda t}A^t+T_t$.  Here $\{A^t\}$ is a $1$-parameter subgroup of $ \mbox{O}(p,q)$, and $\lambda \in \R$.

Let us  begin with a first basic remark.  If $h=e^{\lambda}A+T$ is the time $1$ of the flow $\{h^t\}$, and if $h$ has a fixed point $z_0$ on $\R^{p,q}$, then $\{h^t\}$ has also a fixed point on $\R^{p,q}$.  Indeed, the orbit $h^t.z_0$ is compact.  Its affine convex hull is  compact as well. The flow $\{h^t\}$ is affine and leaves invariant a compact convex set: it must have a fixed point.
We will thus assume in the following that $h=e^{\lambda}A+T$ does not fix any point on $\R^{p,q}$.

If $u \in \R^{p,q}$ is a nonzero lightlike vector, we can associate to $u$ a lightlike geodesic $\beta_u$ through $o$, given by the parametrization:
$$\beta_u(s):=j^o(su)$$
We see $e^{\lambda}A$, and the translation of vector $T$, as elements of $P$ acting on $\Ein^{p,q}$.  Then we can state:

\begin{lemme}
\label{lem.action}
Let $u \in \R^{p,q}$ be a lightlike vector such that $e^{\lambda}A.u=e^{2 \lambda}u$.  Then:
\begin{enumerate}
\item{For every $s \in \R$, $e^{\lambda}A.\beta_u(s)=\beta_u(s).$}
\item{Let $\langle \ ,\    \rangle_{p,q}$ be the bilinear form associated to $Q^{p,q}$.  If $s \not = \frac{1}{\langle T,u \rangle_{p,q}}$, then the translation of vector $T$ maps $\beta_u(s)$ to $\beta_u(\frac{s}{1-s\langle T,u \rangle_{p,q}})$}
\end{enumerate}
\end{lemme}
\begin{preuve}
It is a simple computation from the matrix expressions given in section \ref{sec.matriciel}.
\end{preuve}

The second step is to show:

\begin{lemme}
\label{lem.existence.parabolique}
If $h$ does not fix any point of $\R^{p,q}$, there exists a lightlike vector $u \in \R^{p,q}$ satisfying $e^{\lambda}A.u=e^{2 \lambda}u$, and $\langle T,u \rangle_{p,q} =-1$.
\end{lemme}

\begin{preuve}
we begin with a remark. Assume that $A$ preserves a splitting $\R^{p,q}=F \oplus H$, and assume that $1$ is not an eigenvalue for the restriction of $e^{\lambda}A$ to $H$.  Then, if $T_H$ is the component of $T$ on $H$, there exists $\tau \in H$ such that $(e^{\lambda}A-Id).\tau=T_H$.  Then the translation part of the affine map $(Id+\tau) \circ h \circ (Id +\tau)^{-1}$ does not have any component along $H$.  In other words, conjugating $h$ into $P$, we can assume that $T \in F$.

Let us perform a real Jordan  decomposition of $A$, into the algebraic group $\mbox{O}(p,q)$:  $A=A_sA_eA_u$, where $A_s$, $A_e$ and $A_u$ are all in $ \mbox{O}(p,q)$, $A_s$ is $\R$-semisimple, $A_e$ is elliptic, namely $\C$-semisimple with all its eigenvalues of modulus one, and $A_u$ is unipotent.  Moreover, the elements $A_s$, $A_e$ and $A_u$ are pairwise commuting. Let $F:=\text{Ker}(A_e-Id)$.  Because the group generated by $A_e$ is relatively compact in $\mbox{O}(p,q)$, $A_e$ preserves an orthogonal  splitting $F^{\prime} \oplus H^{\prime}$, where $F^{\prime}$ has dimension $p$ and $Q^{p,q}_{|F^{\prime}}$ is Riemannian, and $H^{\prime}$ has dimension $q$ and $-Q^{p,q}_{|H^{\prime}}$ is Riemannian. As a consequence, $F$ is nondegenerate, of type $(p^{\prime},q^{\prime})$.  We denote by $H$ the orthogonal of $F$, relatively to $\langle \ , \ \rangle_{p,q}$. Then   $A_e$ preserves the splitting $\R^{p,q}=F \oplus H$.  Observe that because $A_s$ and $A_u$ are commuting with $A_e$, $A$ also preserves the splitting $\R^{p,q}=F \oplus H$, and $1$ is not an eigenvalue for the restriction of $e^{\lambda}A$ to $H$.  By the previous remark, we can assume  $T \in F$. 

$\bullet$ {\it Let us first handle  the case where $\lambda=0$}.
The transformation $A_{|F}$ must admit $1$ as an eigenvalue, otherwise $h$ would fix a point in $\R^{p,q}$ by the remark made at the begining of the proof.  Let us call $F_1 \subset F$ the eigenspace associated to $1$.  Because $A_{|F}$ is in $\mbox{O}(p^{\prime},q^{\prime})$, it is not hard to check that $\text{Ker}(A-Id)_{|F}$ is the orthogonal  of $\text{Im}(A-Id)_{|F}$ in $F$.  On the other hand, $T \not \in \text{Im}(A-Id)_{|F}$, otherwise $h$ would have a fixed point in $\R^{p,q}$.  We infer that there exists $u \in \text{Ker}(A-Id)_{|F}$ such that $\langle T,u \rangle_{p,q} \not = 0$.  Rescaling $u$ if necessary, the proposition follows in this case.

$\bullet$ {\it We now assume  that $\lambda \not = 0$.}  In this case, $(A_s)_{|F}$ must admit $e^{-\lambda}$ as an eigenvalue, otherwise $h$ would have a fixed point, by the remark made at the begining of the proof.  Let $F_{1}$ be the associated eigenspace of $e^{\lambda}(A_s)_{|F}$.  Since $(A_s)_{|F}$ is $\R$-semisimple, we can write $F=F_1 \oplus F_2 \oplus \ldots \oplus F_m$, where each $F_j$ is an eigenspace for $(A_s)_{|F}$.  Because $(A_u)_{|F}$ commutes with $(A_s)_{|F}$, this splitting is preserved by $(A_u)_{|F}$, hence by $e^{\lambda}A_{|F}$.  By the remark made at the begining of the proof, we can assume that $T \in F_1$.  Let us observe that because  $A_{|F} \in \mbox{O}(p^{\prime},q^{\prime})$,  $e^{\lambda}$ must also be an eigenvalue of $A_{|F}$, and we will assume that the associated eigenspace is $F_2$.  The sum $F_1 \oplus F_2$ is then nondegenerate of type $(p^{\prime \prime},p^{\prime \prime})$, and  the spaces  $F_1$, $F_2$ are  both totally degenerate (namely, $\langle \ , \ \rangle_{p,q}$ restricts to $0$ on them).  Because $(A_u)_{| F_1 \oplus F_2}$ is in $\mbox{O}(p^{\prime \prime},p^{\prime \prime})$, the spaces $\text{Im}(A_u-Id)_{|F_1 \oplus F_2}$ and $\text{Ker}(A_u-Id)_{|F_1 \oplus F_2}$ are orthogonal.  Since $F_1$ and $F_2$ are left invariant by $A_u$, one has:
$$\text{Im}(A_u-Id)_{|F_1 \oplus F_2}=\text{Im}(A_u-Id)_{|F_1} \oplus \text{Im}(A_u-Id)_{|F_2}$$

Because $h$ does not fix any point in $\R^{p,q}$, $T \not \in \text{Im}(e^{\lambda}A-Id)$, and because $T \in F_1$, we actually get $T \not \in \text{Im}(e^{\lambda}A-Id)_{|F_1}$. On the other hand $(A_u)_{|F_1}=e^{\lambda}A_{|F_1}$, so that $T \not \in \text{Im}(A_u-Id)_{|F_1}$, and finally $T \not \in \text{Im}(A_u-Id)_{|F_1 \oplus F_2}$.  As a consequence, $T \not \in \text{Ker}(A_u-Id)_{|F_1 \oplus F_2}^{\perp}$, hence there exists $u \in \text{Ker}(A_u-Id)_{|F_2}$ such that $\langle T,u \rangle_{p,q} \not =0$.  But   $\text{Ker}(A_u-Id)_{|F_2}=\text{Ker}(e^{\lambda}A-e^{2 \lambda}Id)$, so that rescaling $u$ if necessary, we get the desired proposition. \end{preuve}

As a consequence of lemmas \ref{lem.existence.parabolique} and \ref{lem.action}, we get that if $h$ does not fix any point of $\R^{p,q}$, there exists $u \in \R^{p,q}$ such that for $s \geq 0$, $k \in \N$: 

 $$ h^k.\beta_u(s)=\beta_u( \frac{s}{1+ks})$$
 If $\delta>0$ is small enough, there exists $\alpha : \ ]-\delta , \delta[ \to M$, with $\alpha(0)=x_0$ and  ${\mathcal D}_{x_0}^{b_0}(\alpha)(s)=\beta_u(s)$ for  $s \in \ ]-\delta , \delta[$.

We then use the following reparametrizaton lemma:
\begin{lemme} \cite[Proposition 5.3]{nilpotent}
\label{lem.reparametrage}
Let $I$ be an interval of $\R$ containing $0$, and  $\alpha : I \to M$ a $C^1$-curve  such that $\alpha(0)=x_0$. We call $\beta:={\mathcal D}_{x_0}^{b_0}(\alpha)$.  We assume there exists for every $k \in \N$ a $C^1$-curve  $f_k: I \to I$ satisfying  $f_k(0)=0$, and  such that $h^k.\beta(s)=\beta(f_k(s))$, for $s \in I$.  Then, $\phi_X^k.\alpha(s)=\alpha(f_k(s))$ for $s \in I.$   
\end{lemme}
 
\begin{preuve}
we give an upshot of the proof for the reader's convenience. Let $\hat{\alpha}$ be a lift of $\alpha$ such that $\hat{\alpha}(0)=b_0$, and let $\hat{\beta}:={\mathcal D}_{b_0}(\hat{\alpha})$.  By hypothesis, there exists a $C^1$-curve $p_k:I \to P$, such that $(\Ad h^k).\hat{\beta}(s).p_k(s)=\hat{\beta}(f_k(s))$ for every $s \in I$.  On the other hand, ${\mathcal D}_{b_0}(\hat{\alpha} \circ f_k)= \hat{\beta} \circ f_k$ and 
${\mathcal D}_{b_0}(\phi_X^k.\hat{\alpha}.h^{-k}.p_k)=(\Ad h^k).\hat{\beta}.p_k=\hat{\beta} \circ f_k$.  The two curves $\hat{\alpha} \circ f_k$ and $\phi_X^k.\hat{\alpha}.h^{-k}.p_k$ satisfy the same ODE with the same initial condition: they are equal. \end{preuve}

 Using lemma \ref{lem.reparametrage}, we get that for  $s \in \ [0,\delta[$ and  $k \in \N$, $\phi_X^k.\alpha(s)$ is well defined, and moreover:
 $$ \phi_X^k.\alpha(s)=\alpha(\frac{s}{1+ks})$$
 We are in the situation where $\lim_{k \to \infty}\phi_X^k.\alpha(s)=x_0$ for every $s \in [0,\delta[$, and $D_{x_0}\phi_X^k(\alpha^{\prime}(0))=\alpha^{\prime}(0)$.  Moreover $\alpha^{\prime}(0) \not =0$, because $\beta_u^{\prime}(0) \not = 0$.   The following lemma shows that this forbids  $\{\phi_X^t\}$ being linearizable in a neighborhood of  $x_0$.

\begin{lemme}
\label{lem.obstruction.lineaire}
Let $\{ A^t \}$ be a linear flow on $\R^n$.  Assume that $s \mapsto \alpha(s)$, $s \in ]-\delta, \delta [$ is a $C^1$-curve such that $A^k.\alpha^{\prime}(0)=\alpha^{\prime}(0)$ for every $k \in \N$, and $lim_{k \to +\infty}A^k.\alpha(s)=0$ for every $s \in [0,\delta[$.  Then $\alpha^{\prime}(0)=0$.  

\end{lemme}
\begin{preuve}
we see $\{A^t\}$ as a linear flow on $\C^n$, and call $e^{t\lambda_1}, \ldots , e^{t\lambda_r}$ the eigenvalues of $A^t$, with $\Re e(\lambda_1) \leq \ldots \leq \Re e(\lambda_r)$.  Let $m$ be the greatest integer between $1$ and $r$ such that $\Re e(\lambda_j) <0$ whenever $j \leq m$.  For $i=1, \ldots ,r$, we call $C_i$ the characteristic subspace associated to $e^{t\lambda_i}$, namely $C_i:=\text{Ker}(A^t-e^{t \lambda_i}Id)^n$.  Those spaces do not depend on $t$, and  $\C^n=C_1 \oplus \ldots \oplus C_s$.  This yields a decomposition $\alpha(s)=\alpha_1(s)+ \ldots \alpha_r(s)$, and the condition $\lim_{k \to \infty}A^k.\alpha(s)=0$ for $s \in [0,\delta[$ implies that $\alpha_{m+1}(s)=\ldots=\alpha_r(s)=0$ for $s \in [0,\delta[$.  In particular,  $\alpha(s) \in C_1 \oplus \ldots \oplus C_m$ for $s \in [0,\delta[$, and $\alpha^{\prime}(0) \in C_1 \oplus \ldots \oplus C_m$ as well.  We infer that $\lim_{k \to \infty}A^k.\alpha^{\prime}(0)=0$, and the relation $A^k.\alpha^{\prime}(0)=\alpha^{\prime}(0)$ for all $k \in \N$ yields finally $\alpha^{\prime}(0)=0$.
\end{preuve} 

\begin{remarque}
\label{rem.regularite}
The previous proof shows that if $X$ is smooth, and $\{h^t\}$ does not fix any point in $\R^{p,q}$, then  $X$ is not smoothly linearizable. actually, it shows a little bit more: the field $X$ is  not  even $C^1$-linearizable around $x_0$.
\end{remarque} \end{preuve}

\subsection{Holonomy and essentiality}
We now characterize the local essentiality of $X$, thanks to the holonomy.  Observe that in  \cite[Theorem 2.1]{capocci1}, M. S Capocci obtains the same kind of result, using also the normal Cartan bundle.
\begin{proposition}
\label{prop.essentiel}
The field  $X$ is inessential in a neighborhood of  $x_0$ if and only if its holonomy $\{h^t\} \subset (\R_+^* \times \mbox{O}(p,q)) \ltimes \R^{p,q}$ is conjugated in $P$  to a $1$-parameter group of  ${\mbox O}(p,q)$.

\end{proposition}

 \begin{preuve}
let us assume that  $\{h^t\}$ is conjugated in $P$ to a $1$-parameter group of  ${\mbox O}(p,q)$.  Then, replacing $b_0$ by some  $b_0.p$, we can assume that $\{h^t \} \subset \mbox{O}(p,q)$.  We choose ${\mathcal V}$  a neighborhood of  $0$ in $\lien^-$ such that $\xi \mapsto \pi (\exp(b_0, \xi ))$ is a  diffeomorphism from ${\mathcal V}$ on a neighborhood  $U$ of $x_0$.  Let us call $\Sigma:=\exp(b_0,{\mathcal V})$.  We have a section $\sigma : U \to \Sigma$
 defined by  $\sigma(\pi(\exp(b_0,\xi))):=\exp(b_0,\xi)$, $\forall \xi \in {\mathcal V}$.  
Observe that if  $x=\exp(b_0,\xi)$ is a point of  $ U$, $\xi \in {\mathcal V}$, then for  $t \in \ ]-\delta, \delta[$ small enough, we have: 
$$ \phi_X^t.\sigma(x).h^{-t}=\phi_X^t.\exp(b_0,\xi).h^{-t}=\exp(b_0,(\Ad h^t).\xi)$$
If we choosed  $\delta$ small enough, $\exp(b_0,(\Ad h^t).\xi)$ is a point of  $\Sigma$ projecting on  $\phi_X^t.x$. Thus, it must be  $\sigma(\phi_X^t.x)$ and we get the  equivariance relation:
\begin{equation}
\label{equ.equivariance-section}
 \phi_X^t.\sigma(x).h^{-t}=\sigma(\phi_X^t.x).
\end{equation}
Now, let us fix $\lambda$ a type-$(p,q)$ scalar product on  $\lieg/\liep$, invariant for the adjoint action of ${\mbox O}(p,q)$.  We define on $U$ a type-$(p,q)$ metric  $\mu$ by the following formula:
$$ \mu_x(\iota_{\sigma(x)}(\zeta), \iota_{\sigma(x)}(\zeta)):=\lambda(\zeta,\zeta)$$
Here  $\iota_{\sigma(x)} : \lieg/\liep \to T_xM$ is the isomorphism introduced in section \ref{sec.principe-equivalence}.
By the construction of the normal Cartan bundle (see section \ref{sec.principe-equivalence}), this metric  $\mu$ is in the conformal class $[g]_{|U}$. We are going to show it is invariant by the local flow $\{\phi_X^t\}$, which will prove that  $X$ is a  Killing field for  $\mu$.
Let us recall the  equivariance relation:
\begin{equation}
 \label{equ.equivariance2}
 \iota_{b.p}(\zeta)=\iota_b((\Ad p).\zeta)
 \end{equation}
for $b \in B$, which implies:
\begin{equation}
 \label{equ.differentielle}
 D_x\phi_X^t(\iota_{b}(\zeta))=\iota_{\phi_X^t.b}(\zeta)
 \end{equation}
for $b \in B$ in the fiber of $x$.
We can compute, using  (\ref{equ.equivariance-section}), (\ref{equ.equivariance2}) and (\ref{equ.differentielle}):
$$ \mu_{\phi_X^t.x}(D_x\phi_X^t(\iota_{\sigma(x)}(\zeta)), D_x\phi_X^t(\iota_{\sigma(x)}(\zeta)))=\mu_{\phi_X^t.x}(\iota_{\phi_X^t.\sigma(x)}(\zeta), \iota_{\phi_X^t.\sigma(x)}(\zeta)   ) $$
$$ = \mu_{\phi_X^t.x}(\iota_{\phi_X^t.\sigma(x).h^{-t}}((\Ad h^t).\zeta), \iota_{\phi_X^t.\sigma(x).h^{-t}}( (\Ad h^t).\zeta)   ) $$
 $$ = \mu_{\phi_X^t.x}(\iota_{\sigma(\phi_X^t.x)}((\Ad h^t).\zeta), \iota_{\sigma(\phi_X^t.x)}((\Ad h^t).\zeta))$$
$$ = \lambda((\Ad h^t).\zeta,(\Ad h^t).\zeta))=\lambda(\zeta, \zeta)=\mu_x(\iota_{\sigma(x)}(\zeta),\iota_{\sigma(x)}(\zeta))$$

Reciprocally, if  $X$ is inessential on a neighborhood  $U$ of $x_0$, {\it i.e} $X$ is a  Killing field for a metric  $\mu \in [g]_{|U}$, then $X$  is linearizable at $x_0$.  By the proposition  \ref{prop.lineaire}, $\{h^t\}$ must fix a point of  $\R^{p,q}$.  This means there exists  $p \in P$ such that $p.h^t.p^{-1} \in \R_+^* \times \mbox{O}(p,q)$ for all $t \in \R$.  Thus, replacing $b_0$ by $b_0.p$, we can assume that  $\{ h^t \} \subset \R_+^* \times \mbox{O}(p,q)$. Now, we define on  $\lieg/\liep$ a type-$(p,q)$ scalar product  $\nu$ by the formula:
$$ \nu(\zeta,\zeta):=\mu_{x_0}(\iota_{b_0}(\zeta),\iota_{b_0}(\zeta)) , \ \zeta \in \lieg/\liep$$
Using again (\ref{equ.equivariance-section}), (\ref{equ.equivariance2}) and (\ref{equ.differentielle}), we get that  for every $\zeta \in \lieg/\liep$:
$$ \nu((\Ad h^t).\zeta,(\Ad h^t).\zeta)=\mu_{x_0}(\iota_{b_0}((\Ad h^t).\zeta),(\Ad h^t).\zeta))$$
$$ = \mu_{x_0}(D_{x_0}\phi_X^t(\iota_{b_0}(\zeta)),D_{x_0}\phi_X^t(\iota_{b_0}(\zeta)) )$$
$$ = \mu_{x_0}(\iota_{b_0}(\zeta), \iota_{b_0}(\zeta))=\nu(\zeta,\zeta)$$
We infer that $\nu$ is $(\Ad h^t)$-invariant, what implies that $(\Ad h^t)$ has determinant $1$. The projection of  $\{h^t\}$ on the $\R_+^*$-factor in   $\R_+^* \times \mbox{O}(p,q)$ must then be trivial, so that $\{h^t\} \subset  \mbox{O}(p,q)$. \end{preuve}

\section{Local conformal dynamics of stable sequences}
\label{sec.dynamique-holonomy}
In this section, $(M,[g])$ denotes a smooth, type-$(p,q)$ conformal pseudo-Riemanniann structure of  dimension $n \geq 3$. We call $(M,B,\omega)$ the corresponding normal  Cartan bundle. 

\subsection{Holonomy of a sequence of conformal embeddings}
\label{sec.holonomy-sequence}
The material below is essentially borrowed from \cite[section 4]{familles.normales}.  More details can be found in this reference.

  Let $U$ be an open subset of $M$ and  $f_k : U \to M$ a sequence of conformal embeddings.  Let $x \in U$  be a point, and we assume  that $f_k(x)$ is relatively compact  in $M$.  One says that a sequence   $(h_k)$ of $P$ is a  {\it  holonomy sequence of  $(f_k)$ at  $x$} if there exists a sequence  $(b_k)$ in the fiber of  $x$, contained in a compact set of  $B$, and such that  $f_k(b_k).h_k^{-1}$ is also contained in a compact set of  $B$.

Let us begin with some remarks. First of all, there always exists a holonomy sequence associated to  $(f_k)$.  Also, the notion of holonomy sequence is stable by  ``compact perturbation": if $(h_k)$  is a holonomy sequence of $(f_k)$ at $x$, then so is any sequence $h_k^{\prime}=l_1(k)h_kl_2(k)$, where $l_1(k)$ and $l_2(k)$ are relatively compact sequences of $P$.  One then says that $(h_k)$ and  $(h_k^{\prime})$ are {\it equivalent}.  
Since the action of $P$ on $B$ is proper, it is easily checked that reciprocally, any two holonomy sequences of  $(f_k)$ at $x$  are always equivalent.
Thus, what is really meaningfull is the equivalence class of holonomy sequences of  $(f_k)$ at a point $x$. In what follows, when we will say  {``\it let $(h_k)$ be the holonomy of $(f_k)$ at $x$"}, we will mean that  $(h_k)$ is a representative of the equivalence class of holonomy sequences at $x$. For convenience, we will often be led to change a holonomy sequence into another equivalent one.

\subsection{Notion of stability.}
Let $U$ be an open subset of  $M$ and  $f_k : U \to M$ 
a sequence of conformal embeddings.   
\begin{definition}[Stability]
\label{defi.stabilite}
One says that the sequence  $(f_k)$ is stable at $x \in U$ if for every sequence  $(x_k)$ of $U$ converging to $x$, $f_k(x_k)$ has a same limit  $x_{\infty} \in M$.  The sequence  $(f_k)$ is said to be strongly stable at  $x$ if there exists a neighborhood  $V \subset U$ containing  $x$ such that $f_k(\overline{V})$ converges to $x_{\infty} \in M$ for the Hausdorff topology.
\end{definition}

We will also introduce a notion of stability for sequences in  $P$:
\begin{definition}
\label{defi.opq-stabilite} 
A sequence  $(h_k)$ of $P$ is said to be stable if it is a sequence of $A^+$. Equivalently, $(h_k)$ is stable if it can be written:
$$ h_k=\text{diag}(\lambda_1(k), \ldots , \lambda_n(k)) \in \R_+^* \times \mbox{O}(p,q),$$
with $\lambda_1(k) \geq \ldots \geq \lambda_n(k) \geq 1$.  The sequence $(h_k)$  is said to be strongly stable if moreover all the sequences $\frac{1}{\lambda_i(k)}$ tend to $0$.
\end{definition}

The following lemma shows that the stability at $x$ of  a sequence  $(f_k)$ as above is encoded in its holonomy at $x$. In particular, a sequence $(h_k)$ of $P$, that we see as a sequence of conformal transformations of  $\Ein^{p,q}$, is stable at $o$ (in the sense of definition \ref{defi.stabilite}) if and only if it is stable in the sense of definition \ref{defi.opq-stabilite}. There is thus no ambiguity in the terminology. 

\begin{lemme}\cite[ Lemma 4.3]{familles.normales}
\label{lem.holonomy-stable}
The sequence $(f_k)$ is stable (resp. strongly stable) at $x \in U$  if and only if $f_k(x)$ converges to $x_{\infty} \in M$ and if there exists at  $x$ a holonomy sequence  $(h_k)$ which is  stable (resp.  strongly stable).  
\end{lemme}

%


\section{Semi-completeness properties for conformal vector fields}
\label{sec.semi-complet}
We are still considering a smooth, type-$(p,q)$ conformal structure  $(M,[g])$,  $p+q \geq 3$. We denote by $(M,B,\omega)$ the corresponding normal Cartan bundle. We assume that there exists on  $M$ a conformal vector field   $X$, having a singularity $x_0 \in M$.  We choose $b_0 \in B$ above  $x_0$ and call  $\{h^t\}$ the holonomy flow of $X$ at $x_0$ relatively to  $b_0$.  We are going to use the work done so far to exhibit conditions on  $\{h^t\}$ ensuring that some integral curves of  $X$ in  $M$ are {\it semi-complete}, {\it i.e} defined on $[0,\infty[$ or $]-\infty,0]$.

\subsection{First semi-completeness properties}
\begin{proposition}
\label{prop.semi-complet}
There  exists a real $R_0>0$ such that if $\alpha : [0,1] \to M$ is a conformal geodesic emanating from   $x_0$,  and if $\beta := {\mathcal D}_{x_0}^{b_0}(\alpha)$ satisfies $h^t.[\beta] \subset B(o,R_0)$ for every $t \geq 0$, then $\phi_X^t.\alpha(u)$ is defined  for every $t \geq 0$ and  $u \in [0,1]$. 

\end{proposition}

\begin{preuve}
we first fix $W$  a relatively compact neighborhood  of $x_0$ in $M$,  and we denote by $r_{W}$ the real number given by lemma  \ref{lem.lemme-court}. Then we choose $R_0>0$ such that $8nR_0<r_W$.     Proposition \ref{prop.segment-court} and lemma \ref{lem.lemme-court} ensure that if  $[\gamma]$  is a conformal geodesic segment emanating from $x_0$, and if  ${\mathcal D}_{x_0}^{b_0}([\gamma])=[\zeta]$ is included in  $B(o,R_0)$, then $[\gamma] \subset W$.  

Let us now take a geodesic segment  $[\alpha]$ emanating from  $x_0$, such that $[\beta]:={\mathcal D}_{x_0}^{b_0}([\alpha])$ satisfies $h^t.[\beta] \subset B(o,R_0)$ for every $t \geq 0$. In particular $[\alpha] \subset W$.  We are going to show that if the flow  $\{\phi_X^t\}$ is defined on some interval $[0,T[$, $T>0$ at each  point of $\alpha$, then $\phi_X^t.[\alpha] \subset W$ for every $t \in [0,T[$.  Because $W$ is relatively compact, this will prove that  $\phi_X^t$ is defined  for every $t \geq 0$  at each point of $[\alpha]$.
Since $[\beta]={\mathcal D}_{x_0}^{b_0}([\alpha])$ and  $\{h^t\}$ is the holonomy flow of  $X$ at $b_0$, we have  for every $t \in [0,T[$, ${\mathcal D}_{x_0}^{b_0}(\phi_X^t.[\alpha])=h^t.[\beta]$.
But by assumption $h^t.[\beta] \subset B(o,R_0)$ for $t \geq 0$, and in particular for $t \in [0,T[$ so that we get  $\phi_X^t.[\alpha] \subset W$ for every $t \in [0,T[$,  concluding the proof.
\end{preuve}

From this proposition, we infer:

\begin{corollaire}
\label{coro.compact-complet}
Let $\liei_{x_0}$ be the Lie algebra of conformal vector fields on $M$ vanishing at $x_0$, and let $\lieh \subset \liei_{x_0}$ be a subalgebra.  Let $H_h \subset P$ be the holonomy group of $\lieh$ with respect to $b_0$. If $H_h$ is relatively  compact in $P$,  then there exists a neighborhood  $U$ of $x_0$  on which the local action of $\lieh$ integrates into the action of a connected Lie subgroup $H \subset \Conf(U)$, which is relatively compact in $\Conf(U)$.  

\end{corollaire}

\begin{preuve}
observe first that any compact subgroup of $P$ is conjugated in $P$ to a subgroup of $\mbox{O}(p,q)$.  Thus, replacing $b_0$ by a suitable $b_0.p$, we will assume that $(\Ad H_h)$ leaves $\lien^-$ invariant.  We  define ${\mathcal U}_r:=\bigcup_{h \in H_h} (\Ad h).{\mathcal B}(0,r)$.  This is an open subset of  $\lien^-$ and because  $H_h$ is relatively compact, $\lim_{r \to 0}{\mathcal U}_r=\{0\}$ (the limit being taken for the Hausdorff topology).  Let  $ r>0$  be small enough so that on the one hand $\xi \mapsto \pi(\exp(b_0,\xi))$ is a diffeomorphism from  ${\mathcal U}_r$ onto its  image  $U$,  and on the other hand   $\pi_G(\exp({{\mathcal U}_r}) ) \subset B(o,R_0)$, where $R_0>0$ is given by  proposition \ref{prop.semi-complet}. 

We first show  that the local action of $\lieh$ integrates into the action of a Lie group on $U$.  This amounts to show that any vector field $X \in \lieh$ is complete on $U$.  Let $X \in \lieh$ and $\{h^t\}$ its holonomy flow with respect to $b_0$.
   We have, for every $u \in [0,1]$,  and every  $\xi \in {\mathcal U}_r$:
$$ {\mathcal D}_{x_0}^{b_0}(\alpha)(u)=\pi_G(\exp({u\xi})),$$
where  $\alpha(u):=\pi(\exp(b_0,u \xi))$.  Also,  $h^t.\pi_G(\exp({u \xi}))=\pi_G(\exp({u(\Ad h^t)).\xi}).$  Since $(\Ad h^t).{\mathcal U}_r= {\mathcal U}_r$, we obtain that  $h^t.\pi_G(\exp({u \xi})) \in B(o,R_0)$ for all $u \in [0,1]$, $t \geq 0$.  Proposition \ref{prop.semi-complet} then ensures that for  every $x \in U$, $\phi_X^t.x$ is defined  for every $t \geq 0$.  Considering  $-X$ instead of  $X$, we obtain in the same way that   $\phi_X^t.x$ is defined  for every  $t \leq 0$.  

Let us call $H \subset \Conf(U)$ the connected Lie subgroup having $\lieh_{|U}$ as Lie algebra. We are going to show that $H$ is relatively compact in $\Conf(U)$, for the topology of uniform convergence on compact subsets of $U$. Let $(\phi_k)$  be a sequence of $H$.  There exists an holonomy sequence $(h_k) \subset H$ such that $\phi_k.b_0.h_k^{-1}=b_0$.  By assumption on $H$, there exists a subsequence $(h_{k_j})$ which converges to $l \in P$.

The  relation:
$$ \phi_{k_j}.\exp(b_0,\xi).h_{k_j}=\exp{(b_0,(\Ad h_{k_j}).\xi})$$
shows that  $(\phi_{k_j})$ converges uniformly  on  the compact subsets of $U$ to the diffeomorphism:
$$ \pi(\exp(b_0,\xi)) \mapsto \pi(\exp(b_0. (\Ad l).\xi)), \ \xi \in {\mathcal U}_r$$
(actually the convergence is $C^{\infty}$). We conclude that  $H$ is relatively compact  in the conformal group of  $U$.
\end{preuve}

\subsection{Stability and semi-completeness for conformal vector fields.}
\label{sec.stabilite-completude}
Proposition \ref{prop.semi-complet} shows how some  properties of the holonomy give rise to semi-complete orbits for a conformal vector field.  The good point with the notion of stability introduced in section  \ref{sec.dynamique-holonomy} is that starting from a semi-completeness property for {\it one} orbit, it yields a semi-completeness property on a non empty open set of orbits. One can then study the dynamical behavior  of  $\{\phi_X^t\}$ on this open set, and try to infer some geometrical consequences.  We keep in the statement below the notations introduced so far. 
\begin{proposition}
\label{prop.fort-stable-complet}
 Let $\alpha: [0,1] \to M$ be a conformal geodesic emanating from  $x_0$, $x:=\alpha(1)$, and  $\beta:={\mathcal D}_{x_0}^{b_0}(\alpha)$.  We assume that $\phi_X^t.\alpha(u)$ is defined for every $(t,u) \in \R_+ \times [0,1]$, that  $\lim_{t \to \infty} h^t.[\beta]=o$, and that  for every sequence $t_k \to \infty$,  $(h^{t_k})$  is stable (resp. strongly  stable) at  $\beta(1)$. Then 
\begin{enumerate}
\item{There exists an open set $V$ containing  $x$ such that for every $y \in V$, $\phi_X^t.y$ is defined for every $t \geq 0$.}
\item{Moreover there exists a sequence $(s_k)$ of $\R_+$ tending to infinity, such that for every $y \in V$, $\lim_{k \to \infty}\phi_X^{s_k}.y$ exists, $(\phi_X^{s_k})$ is  stable at $y$  (resp.  $\phi_X^{s_k}.V \to x_0$, and  $(\phi_X^{s_k})$ is strongly  stable at each  point of $V$), and the holonomy of $(\phi_X^{s_k})$  at  $y$ is that of $(h^{s_k})$ at  $\beta(1)$.}
\end{enumerate}

\end{proposition}

\begin{preuve}
by corollary  \ref{coro.tend-vers-0}, we have $\lim_{t \to \infty} \phi_X^t. \alpha(u)=x_0$ for every $u \in [0,1]$. Let us fix a compact set  $K$ of $M$ containing  $\{  \phi_X^t.x \}_{t \geq 0} \cup \{  x_{0}\}$ in its interior.  We want to show that there exists an open neighborhood $V$ of  $x$ having the property that for  every $t_0 >0$ such that $\phi_X^t$ is defined on  $[0,t_0[$  at each point of $V$, we have actually  $\phi_X^t(V) \subset K$, $\forall t \in [0,t_0[$.  

If it is not the case, we can find a sequence $(y_k)$ tending to  $x$, as well as times  $t_k \in [0,\infty[$ such that $\phi_X^t.y_k$ is defined on $[0,t_k]$, but  $\phi_X^{t_k}.y_k \not \in K$.  We are now going to use the following lemma, the proof of which we  postpone a little bit later:

\begin{lemme}
\label{lem.extension-holonomy}
  Let us assume there exists a conformal geodesic  $\alpha: [0,1] \to M$ emanating from  $x_0$ such that for every $u \in [0,1]$, $\phi_X^t.\alpha(u)$ is defined for every $t \geq 0$ (resp. $t \leq 0$), and  $\lim_{t \to \infty} h^t.[\beta]=o$ (resp. $\lim_{t \to -\infty} h^t.[\beta]=o$), where $\beta:={\mathcal D}_{x_0}^{b_0}(\alpha)$.  Then  for every sequence $t_k \to \infty$  (resp. $t_k \to -\infty$), and every $u \in [0,1]$, the holonomy of  $(\phi_X^{t_k})$ at  $\alpha(u)$ is the holonomy of  $(h^{t_k})$ at  $\beta(u)$. 
\end{lemme}

This lemma tells us that the holonomy  $(h_k)$  of $(h^{t_k})$ at  $\beta(1)$ is also an holonomy sequence of  $(\phi_X^{t_k})$ at $x$.  Because  $(h^{t_k})$ is stable at $\beta(1)$, we can assume that $(h_k)$ is a sequence of $A^+$ that we write $h_k=\diag(\lambda_1(k), \ldots, \lambda_n(k))$.  Let $(b_k)$ be a relatively compact  sequence of  $B$, such that $\phi_X^{t_k}(b_k).h_k^{-1}$ is also relatively compact in    $B$.  Because $(y_k)$ tends to $x$, there is a  sequence $(\xi_k)$ of $\lien^-$ tending to  $0$, and satisfying  $\pi(\exp(b_k,\xi_k))=y_k$.   We have:
$$ \phi_X^{t_k}(\exp(b_k,\xi_k).h_k^{-1})=\exp(\phi_X^{t_k}(b_k).h_k^{-1},(\Ad h_k).\xi_k)$$ 
The stability of $(h^{t_k})$ at  $\beta(1)$  means that $|\frac{1}{\lambda_i(k)}|$, $i=1, \ldots, n$ are bounded sequences, ensuring that   $(\Ad h_k).\xi_k$ tends to  $0$.  Projecting on $M$, we get that $\phi_X^{t_k}(y_k)$ tends to $x_{0}$.   But  $x_{0}$ is in the interior of  $K$, contradicting  $\phi_X^{t_k}(y_k) \not \in K$.  This shows the first  point of the  proposition. 

Let us now prove the  second point.  We first choose a sequence $(t_k)$ tending to infinity.  By the first point, we know that $\phi_X^{t_k}$ is well defined on a neighborhood of $x$ for all $k \in \N$, and  lemma \ref{lem.extension-holonomy} ensures that $(\phi_X^{t_k})$ is stable (resp. strongly stable) at $x$.  In particular, it is equicontinuous and we can apply the theorem 1.1 of \cite{familles.normales}:   replacing if necessary $V$ by one of its open subsets, there exists a sequence $(s_k)$ tending to infinity, which is a subsequence of $(t_k)$, and such that $(\phi_X^{s_k})$ tends uniformly on the compact subsets of $V$ to a smooth map $\phi:V \to M$.  In particular, for every $y \in M$, $\lim_{k \to \infty}\phi_X^{s_k}.y$ exists.  Lemma 6.1 of \cite{familles.normales} ensures that the holonomy $(h_k)$ of $(\phi_X^{s_k})$ at $x$ is an holonomy sequence of $(\phi_X^{s_k})$ at $y$ for every $y \in V$.  Moreover, lemma \ref{lem.extension-holonomy} says that this holonomy $(h_k)$ is that of $(h^{s_k})$ at $\beta(1)$.
\end{preuve}

\subsubsection{Proof of  Lemma  \ref{lem.extension-holonomy}}
We write  $\alpha(u)=\pi(\exp(b_0,u \xi))$. By hypothesis, $u \mapsto \beta_k(u):=h^{t_k}.\pi_G(\exp({u\xi}))$ is a sequence of conformal geodesics tending to   $o$ for the  Hausdorff topology.  Proposition \ref{prop.segment-court} ensures that the length $L^o(\beta_k)$ of those curves tends to  $0$ as $k \to \infty$.  For each  $k \in \N$, there is a $C^1$-curve  $p_k : [0,1] \to P$, with $p_k(0)= 1_G$, such that ${\hat \beta}_k(u):=\exp({u (\Ad h^{t_k}).\xi}).p_k(u)^{-1}$ is a curve of  $N^-$.  As already noticed, the length of   ${\hat \beta}_k$ with respect to $\rho^G$ is just  $L^o(\beta_k)$. We thus get  
$$\lim_{k \to \infty}h^{t_k}.\exp({u \xi}).(p_k(u)h^{t_k})^{-1}=1_G$$
 proving that  $(p_k(u)h^{t_k})$ is a holonomy  sequence of  $(h^{t_k})$ at $\beta(u)=\pi_G(\exp({u\xi}))$.  On the other hand, if we define ${\hat \alpha}_k(u):=\phi_X^{t_k}.\exp(b_0,u\xi).(p_k(u)h^{t_k})^{-1}$, we also have  ${\hat \beta}_k={\mathcal D}_{x_0}^{b_0}({\hat \alpha}_k)$.  Thus the length of ${\hat \alpha}_k$ relatively to  $\rho^B$ tends to $0$ because it is the length of ${\hat \beta}_k$ relatively to $\rho^G$.  Hence  
 $$\lim_{k \to \infty}\phi_X^{t_k}.\exp(b_0,u\xi).(p_k(u)h^{t_k})^{-1}=b_0, \text{ for every } u \in [0,1]$$  This shows that  $(p_k(u)h^{t_k})$ is also a holonomy sequence for  $(\phi_X^{t_k})$ at  $x_0$.

\section{Application to Riemannian conformal vector fields: \\
proof of theorem \ref{thm.ferrand-champs}}
\label{sec.ferrand-champs}

Here, $(M,[g])$ is a smooth Riemannian conformal structure of  dimension $\geq 3$.  We still denote by  $(M,B,\omega)$ the corresponding normal  Cartan bundle.  Let us point out that in the Riemannian framework, the space  $\Ein^{0,n}$  is nothing but the standard conformal sphere ${\bf S}^n$, and  $\R^{0,n}$ is the Euclidean space $\R^n$.
We consider $x_0 \in M$, and  ${\mathfrak I}_{x_0}$ the Lie algebra of smooth conformal vector fields on $M$ vanishing at $x_0$.  We will call ${\mathfrak I}_h$ the associated holonomy algebra with respect to some point $b_0 \in B$ in the fiber of $x_0$, and $I_h \subset P$ will denote  the connected Lie subgroup associated to $\liei_h$.

The aim of this section is to prove the theorem below, of which theorem \ref{thm.ferrand-champs} is a particular case:

\begin{theoreme}
\label{thm.ferrand-algebres}
Let $(M,g)$ be a smooth Riemannian manifold of dimension  $n \geq 3$, and $x_0 \in M$. Let ${\mathfrak I}_{x_0}$ be the Lie algebra of conformal vector fields vanishing at $x_0$, and $\lieh \subset \liei_{x_0}$ a subalgebra. We denote by $\lieh_h$ the holonomy algebra of $\lieh$. Then :

\begin{enumerate}
\item{Either there is a neighborhood  $U$ of $x_0$ on which  the local action of $\lieh$ integrates into the action of a  relatively compact subgroup $H \subset \Conf(U)$.   In this case  $H$ preserves a metric in the conformal class $[g]_{|U}$, and its action around $x_0$ is linearizable.}
\item{If we are not in the previous case, there is a neighborhood  $U$ of $x_0$ which is conformally flat , and the Lie algebra $\lieh$ is essential on each neighborhood of $x_0$.}
\end{enumerate}
In any case, there is a smooth diffeomorphism from a neighborhood of $x_0$ in $M$ onto a neighborhood of $o$ in ${\bf S}^n$, conjugating $\lieh$ and  its holonomy algebra $\lieh_h$. 
\end{theoreme} 

Before going further, let us precise that a subalgebra $\lieh \subset \liei_{x_0} $ is said to be inessential on some neighborhood $U$ of $x_0$ if  there exists a metric $g$ in the conformal class $[g]_{|U}$ such that every vector field of $\lieh$ is a Killing field for $g$ ({\it i.e} generates a local flow of isometries of $g$).  If we are not in this case, the Lie algebra $\lieh$ is said essential.

\subsection{Strengthening of some results in the Riemannian setting}
We begin by proving some technical results allowing to sharpen some of our previous statements in the Riemannian case.

\begin{lemme}
\label{lem.algebrique}
Let $m \geq 1$, and $(\R \oplus \oo (m)) \ltimes \R^m$ be the Lie algebra of  $(\R_+^* \times \mbox{O}(m)) \ltimes \R^m$, the group of affine conformal transformations of $\R^m$.  Let $\lieh$ be a subalgebra of $(\R \oplus \oo (m)) \ltimes \R^m$, and $H$ the connected Lie subgroup of $(\R_+^* \times \mbox{O}(m)) \ltimes \R^m$ having $\lieh$ as Lie algebra.  
\begin{enumerate}
\item{Assume that   for every $Y \in {\lieh}$, the $1$-parameter subgroup $\{ \exp(tY)\}_{t \in \R}$ is relatively compact in $(\R_+^* \times \mbox{O}(m)) \ltimes \R^m$.  Then $H$ is relatively compact in $(\R_+^* \times \mbox{O}(m)) \ltimes \R^m$. }
\item{Assume that every element of $H$ has a fixed point in $\R^m$. Then $H$ has a  fixed point in $\R^m$.}
\end{enumerate}
\end{lemme}

\begin{preuve}
we begin with the first point, that we  prove by induction.  When $m=1$, the claim is obvious.  Assume now $m>1$, and take $Y \in \lieh$.    Because the adjoint representation of $\{\exp(tY)\}$ on $(\R \oplus \oo (m)) \ltimes \R^m$ is relatively compact, the eigenvalues of $(\add Y)$ are purely imaginary, and the same is true for the eigenvalues of $(\add Y)_{|\lieh}$.  As a consequence, if the center of $\lieh$ is trivial, then the Killing form of $\lieh$ is negative definite.  It is then a classical result that $H$ is  compact, as a Lie group, hence is a compact subgroup of $(\R_+^* \times \mbox{O}(m)) \ltimes \R^m$.  

If the center $\liez \subset \lieh$ is nontrivial, and $Z \subset (\R_+^* \times \mbox{O}(m)) \ltimes \R^m$ is the associated connected subgroup, then $\overline{Z}$, the closure of $Z$ into $(\R_+^* \times \mbox{O}(m)) \ltimes \R^m$,  is a compact torus.  In particular, it must fix at least one point of $\R^m$, and the affine subspace $\mathcal{F} \subsetneq \R^m$ of fixed points of $ \overline{Z}$ is acted upon conformally by $H$.    Let $m^{\prime}$ be the dimension of $\mathcal {F}$.  If $m^{\prime}=0$, then $H$ fixes a point in $\R^m$, hence is conjugated to a subgroup of $\R_+^* \times \mbox{O}(m)$.  Because all the $1$-parameter subgroups of $H$ are relatively compact, $H$ is actually conjugated to a subgroup of $\mbox{O}(m)$, hence is relatively compact in $(\R_+^* \times \mbox{O}(m)) \ltimes \R^m$.

Now, if $m^{\prime} \geq 1$, we can use the induction hypothesis:  the restriction $H_{|\mathcal{F}}$ is relatively compact into $(\R_+^* \times \mbox{O}(m^{\prime})) \ltimes \R^{m^{\prime}}$.  We infer that $H$ fixes a point of $\mathcal {F}$, and the same argument as above implies that $H$ is conjugated to a subgroup of $\mbox{O}(m)$.

We now prove the second point of the lemma, arguing again by induction.  For $m =1$, the claim is pretty clear.  Assume now that $m>1$.  We consider the Lie algebra homomorphism $\pi_1: \lieh_h \subset (\R \oplus \oo (m)) \ltimes \R^m \to \R$, and call $\lieh_h^{\prime}$ the kernel of $\pi_1$.  If this kernel is trivial, then $H_h$ is either trivial, or is conjugated to  a $1$-parameter subgroup  of $\R_+^* \times \mbox{O}(m)$.  In any case, it fixes a point of $\R^m$.  

We assume now that $\lieh_h^{\prime}$ is nontrivial, and call $H_h^{\prime}$ the connected subgroup of  $ \mbox{O}(m) \ltimes \R^{m}$ with Lie algebra $\lieh_h^{\prime}$.  By hypothesis, if $Y \in \lieh_h^{\prime}$, then $\{ \exp(tY) \}$ has a fixed point in $\R^m$.  But a $1$-parameter subgroup of $ \mbox{O}(m) \ltimes \R^{m}$ having a fixed point is conjugated to a $1$-parameter subgroup of $ \mbox{O}(m)$, hence is relatively compact.  We infer from the first point of the lemma that $H_h^{\prime}$ is relatively compact in $ \mbox{O}(m) \ltimes \R^{m}$, hence fixes a point in $\R^m$.  We call $\mathcal{F} \subsetneq \R^m$ the affine subspace, of dimension $m^{\prime}<m$, comprising the fixed points of $H_h^{\prime}$.  Because $H_h^{\prime}$ is normal in $H_h$,  $\mathcal{F}$ is left invariant by $H_h$, and $(H_h)_{|\mathcal{F}}$ is a subgroup of $(\R_+^* \times\mbox{O}(m^{\prime})) \ltimes \R^{m^{\prime}}$. We are going to show that every $1$-parameter subgroup $\{ \exp(tY ) \}$, $Y \in \lieh_h$, has a fixed point in $\mathcal{F}$.  It is clear if $Y \in \lieh_h^{\prime}$.  Now, if $Y \in \lieh_h \setminus \lieh_h^{\prime}$, then $\{ \exp(tY ) \}$ is conjugated to a subgroup $e^{\lambda t} A^t$, where $\lambda \not =0$, and $\{ A^t \} \subset \mbox{O}(m)$.  In particular, $\{ \exp(tY) \}$ has a unique fixed point $y_0$ in $\R^n$, having the property that $\lim_{t \to \infty} \exp(tY).y=y_0$ for every $y \in \R^m$, or $\lim_{t \to -\infty} \exp(tY).y=y_0$ for every $y \in \R^m$.  Choosing $y \in \mathcal{F}$, and because $\mathcal{F}$ is a closed subset, we get $y_{0} \in \mathcal{F}$.  We can then use the induction hypothesis  and obtain that $H_h$ fixes a point in $\mathcal{F}$.  This concludes the proof of the second point. \end{preuve}

As a corollary of  this lemma, we can strengthen the statements of propositions \ref{prop.lineaire} and \ref{prop.essentiel}:

\begin{corollaire}
\label{coro.lineaire.global}
Let $\lieh \subset \liei_{x_0}$ be a subalgebra of smooth conformal vector fields.  Let $b_0 \in B$ in the fiber of $x_0$ and  $\lieh_h$ (resp. $H_h$)  the  holonomy algebra (resp. the holonomy group) of $\lieh$ with respect to $b_0$.   Then:
\begin{enumerate}
\item{The algebra  $\lieh$ is linearizable in a neighborhood of $x_0$ if and only if $H_h$ has a fixed point in $\R^n$. In this case, there is a smooth diffeomorphism from a neighborhood of $x_0$ in $M$  onto a neighborhood of $o$ in ${\bf S}^n$ which conjugates $\lieh$ and $\lieh_h$.}
\item{The algebra $\lieh$ is inessential in a neighborhood of $x_0$ if and only if $H_h \subset P$ is conjugated into $P$ to a subgroup of $\mbox{O}(n)$.}
\end{enumerate}
\end{corollaire}

\begin{preuve}
we keep the notations of proposition \ref{prop.lineaire}.  If $H_h$  fixes a point of $\R^n$, we may assume that this point is $0$.  Then the smooth diffeomorphism $\psi \circ \varphi^{-1}$, from $U$ onto $V$ conjugates $\lieh$ and $\lieh_h$.

Reciprocally, if $\lieh$ is smoothly linearizable around $x_0$, proposition \ref{prop.lineaire} ensures that for every $Y \in \lieh_h$, the $1$-parameter group $\{ \exp(tY) \}_{t \in \R}$ has a fixed point in $\R^n$.  Then we can apply the second point of lemma \ref{lem.algebrique} and infer that $H_h$ fixes a point in $\R^n$.

We now prove the second point of the corollary.  We can reproduce {\it verbatim} the proof of proposition \ref{prop.essentiel} to get that whenever $H_h$ is conjugated to a subgroup of $\mbox{O}(n)$, then $\lieh$ is inessential in a neighborhood $U$ of $x_0$.   Reciprocally, if $\lieh$ is inessential in a neighborhood of $x_0$, it is smoothly linearizable around $x_0$.  By the first point of the corollary, $H_h$ must fix a point of $\R^n$, hence is conjugated to a subgroup of $\R_+^* \times \mbox{O}(n)$. Because each $X\in \lieh$ is inessential around $x_0$, every $1$-parameter group of $H_h$ has determinant $1$, showing that 
$H_h$ is actually conjugated to a subgroup of $\mbox{O}(n)$. \end{preuve}

We can now prove theorem \ref{thm.ferrand-algebres}  in  the case when $H_h$ is relatively compact in $P$. In this case, $H_h$  has a fixed point in $\R^n$, so that by corollary   \ref{coro.lineaire.global}, $\lieh$ is smoothly linearizable and there exists a smooth diffeomorphism  $\varphi$ from a neighborhood  $U$ of  $x_0$ onto a neigborhood  $V$ of $o$, conjugating $\lieh$  and the holonomy algebra  $\lieh_h$. The second point of corollary  \ref{coro.lineaire.global} says that  $\lieh$ is inessential  in a neighborhood of  $x_0$. Finally, corollary \ref{coro.compact-complet} ensures that there exists a neighborhood $U$ of $x_0$ on which  every vector field of $\lieh$ is complete. In other words,  the local action of $\lieh$ integrates into the action of a Lie subgroup $H \subset \Conf(U)$.  We are then in case  $(1)$ of theorem  \ref{thm.ferrand-algebres}.

Whenever  $H_h$ is not relatively compact in $P$, lemma \ref{lem.algebrique} ensures that for some $X_h \in \lieh_h$, the $1$-parameter group $\{h^t \}:=\{ \exp(t X_h) \}$ is not relatively compact in $P$.  There exists $X \in \lieh$ such that $\{ h^t \}$ is the holonomy flow of $X$ with respect to $b_0$.  Up to conjugacy in $P$, the flow $\{ h^t \}$ is of two possible kinds, and we are going to prove theorem \ref{thm.ferrand-algebres} in each case.

\subsection{The flow  $\{h^t\}$ is a non relatively compact  flow of $\R_+^* \times \mbox{O}(n)$}

The adjoint action of  $h^t$ on  $\lien^-$ is by dilations for the metric  $< \ , \ >_{\lien^-}$.  In particular there exists  $\lambda  \not = 0$ such that $(\Ad h^t).{\mathcal B}(0,r)={\mathcal B}(0,e^{\lambda t}r)$, for every $r>0$.  Replacing if necessary  $X$ by $-X$, we will assume that  $\lambda < 0$.  Proposition \ref{prop.lineaire} ensures that  $X$ is locally conjugated in a neighborhood of $x_0$ to its holonomy field $X_h$, which is a field of linear contractions in $\R^n$.  Thus, there exists a  neighborhood $U$ of  $x_0$ such that  for every point $x \in U$, $\phi_X^t.x$ is defined for every $t \geq 0$, and moreover  $\lim_{t \to \infty} \phi_X^t.U=x_0$.  In particular the  sequence $\{\phi_X^k \}_{k \geq 0}$ is strongly  stable  at each  point of  $U$.

We conclude that  $U$ is conformally flat thanks to the:

\begin{lemme}
\label{lem.conf-plat}
Let $(M,[g])$ be a conformal Riemannian structure of  dimension $n \geq 3$, and  $f_k : U \to M$ a sequence of conformal embedings,  which is strongly stable at  $x \in U$. Then the  Weyl  tensor (resp. the Cotton tensor if  $\text{dim }M=3$)  vanishes on a neighborhood of  $x$.
\end{lemme}
\begin{preuve}
If $f_k:U \to M$ is strongly stable at $x$, then it is equicontinuous at $x$.  We can then use \cite[Proposition 5]{frances-ferrand} (the proof of which deals with conformal diffeomorphisms of $M$, but can be adapted to conformal embeddings in a straigthforward way) and get that the  Weyl  tensor (resp. the Cotton tensor if  $\text{dim }M=3$)  vanishes at $x$.  By definition, being strongly stable is an open property, so that the lemma follows.\end{preuve}

We are thus in case  $(2)$  of theorem \ref{thm.ferrand-algebres}.  Proposition \ref{prop.essentiel} ensures that  $X$, hence $\lieh$, is essential on any neighborhood of  $x_0$. By the remark  \ref{rem.cas-plat}, $\lieh$ is smoothly conjugated (actually the conjugacy map is  a smooth {\it conformal} diffeomorphism) to its holonomy  algebra $\lieh_h$.

\subsection{The flow $\{h^t\}$ is the commutative product of a   translation  flow ${\tau^t}$, and a relatively compact flow  ${\kappa}^t$ in  ${\mbox O}(n)$.}

  For every $R >0$, and $z=\pi_G(\exp({\xi})) \in B(o,R)$, $\xi \in {\mathcal B}(0,R)$, we call  $[oz]$  the geodesic segment defined by:
$$ [oz]:=\{ \pi_G(\exp({u \xi})) \ | \ u \in [0,1]  \}$$

\begin{lemme}

\label{lem.piege}
 For every $R>0$ and $z \in B(o,R)$, we are in one of the following cases:
\begin{itemize}
\item{For every  $t \geq 0$, $h^t.[oz] \subset B(o,R)$ and $\lim_{t \to \infty} h^t.[oz]=o$.}
\item{For every  $t \leq 0$, $h^t.[oz] \subset B(o,R)$ and $\lim_{t \to -\infty} h^t.[oz]=o$.}
\end{itemize}
\end{lemme}

\begin{preuve}
with the notations of  section \ref{sec.metrique-auxiliaire}, and in the Riemannian case, we have $<x,x>=||x||^2=Q^{0,n}(x)$.  Everywhere it is defined, the map  $ s(x)=\frac{-2x}{<x,x>} $  satisfies:
$$j^o(s(x))=j(x)$$

We first prove the proposition when  $h^t$ is the translation of vector  $tv$, for a nonzero $v \in \R^n$.
We write  $z=j(x)$, and  $z \in B(o,R)$ means that we have $||s(x)|| < R$, {\it i.e} $||x||> \frac{2}{R}$.
As a set, the  segment $[oz]$ is:
$$[oz]=\{o \} \cup \{ j(ux) \ | \ u \in [1,\infty[ \}$$

From $|| s(ux+tv) ||^2=4(u^2||x||^2+t^2||v||^2+2ut<x,v>)^{-1}$, we infer:

- if $<x,v> \ \geq 0$, then $|| s(ux+tv) ||^2 \leq \frac{4}{||x||^2}$ for every $t \geq 0$,  meaning that   $h^t.[oz] \subset B(o,R)$ for  $t \geq 0$.  Also, from  $\sup_{u \in [1,\infty[}|| s(ux+tv) ||^2 \leq \frac{1}{t^2||v||^2}$ for $t>0$, we get  $\lim_{t \to \infty} \tau^t.[oz]=o$.

- if $<x,v> \ \leq 0$, then  $|| s(ux+tv) ||^2 \leq \frac{4}{||x||^2}$ for every $t \leq 0$, meaning that  $h^t.[oz] \subset B(o,R)$ for $t \leq 0$.  Also, from  $\sup_{u \in [1,\infty[}|| s(ux+tv) ||^2 \leq \frac{1}{t^2||v||^2}$ for $t<0$, we get  $\lim_{t \to -\infty} \tau^t.[oz]=o$.

If now  $\{h^t\}$  is the commutative product of $\{\tau^t\}$,  the translation flow  of direction $tv$,  and of a flow  $\{ \kappa^t \}$ in ${\mbox O}(n)$, then the lemma remains true since $\kappa^t.[oz]=[o\kappa^t.z]$, and the balls  $B(o,R)$ are left invariant  under the action of  $\kappa^t$.  
\end{preuve}

\begin{lemme}
\label{lem.translation-stable}
Let $(t_k)$ be a   sequence  of real numbers tending to  $\infty$  (resp. to $- \infty$).  Then the sequence $(h^{t_k})$   is strongly stable  at each point $z \in j(\R^n)$.
\end{lemme}

\begin{preuve}
we write that $h^{t_k}$  is a commutative product  $\tau^{t_k}.\kappa^{t_k}$, with $\tau^{t_k}$ the translation of vector $t_kv \not = 0$, and  $\kappa^{t_k} \in \mbox{O}(n)$.
Let $z=j(x) \in j(\R^n)$ and  $y$  be in the euclidean ball of center $x$ and radius  $r>0$.  Then:
$$ \lim_{k \to \infty} || s(y+t_kv)||^2=0$$
uniformly on  $B(x,r)$.   As a consequence  $\lim_{k \to \infty} \tau^{t_k}.j(B(x,r))=o$.  Because  $\{ \kappa^{t_k} \}$ has compact closure in  $P$, we have also:
$$ \lim_{k \to \infty} h^{t_k}.j(B(x,r))=o,$$
showing that  $(h^{t_k})$  is strongly stable at $z$.  
\end{preuve}


Let us consider  $R_0$ the real number given by proposition  \ref{prop.semi-complet}, and let us choose $0<R<R_0$ small enough so that  $\xi \mapsto \pi(\exp(b_0,\xi))$ is defined and  injective on  ${\mathcal B}(0,R)$.   Let us define  $U:=\pi ( \exp(b_0,{\mathcal B}(0,R)))$.  Pick   $x \in U$, and write  $x=\pi (\exp(b_0,\xi))$ for  $\xi \in {\mathcal B}(0,R)$.  Setting   $\alpha(u):=\pi(\exp(b_0,u \xi))$,  we see thanks to  lemma \ref{lem.piege} that  $\beta:={\mathcal D}_{x_0}^{b_0}(\alpha)$ satisfies either $h^t.[\beta] \subset B(o,R_0)$ for every $t \geq 0$, and $\lim_{t \to \infty} h^t.[\beta]=o$, or  $h^t.[\beta] \subset B(o,R_0)$ for every $t \leq 0$, and $\lim_{t \to -\infty} h^t.[\beta]=o$. We will assume in the following that we are in the first of those two cases (the proof is the same for the second case).  Proposition \ref{prop.semi-complet}  ensures that  $\phi_X^t.\alpha(u)$ is defined for every $t \geq 0$  and  $u \in [0,1]$.  Corollary  \ref{coro.tend-vers-0} yields  $\lim_{t \to \infty} \phi_X^t.\alpha(u)=x_0$ for every $u \in [0,1]$.  
We can now use  proposition \ref{prop.fort-stable-complet}: there exists  $V$ a neighborhood of  $x$ in $U$ such that $\phi_X^t$ is defined on  $V$ for every $t \geq 0$.  On the other hand, the  second point  of  proposition  \ref{prop.fort-stable-complet}, together with the third point of  lemma  \ref{lem.translation-stable}   ensures that for some sequence $s_k \to \infty$, $(\phi_X^{s_k})$ is strongly  stable at  $x$.  Lemma \ref{lem.conf-plat} then   says that a neighborhood of $x$ must be conformally flat.  This holds  for every $x \in U$, thus $U$ itself is conformally flat.  Thanks to the remark \ref{rem.cas-plat}, $\lieh$ is locally conjugated by a smooth conformal diffeomorphism to its holonomy algebra $\lieh_h$. Proposition \ref{prop.essentiel} ensures that the field  $X$, hence the algebra $\lieh$,  is essential on any neighborhood of $x_0$. Once again, we are in case  $(2)$ of theorem \ref{thm.ferrand-algebres}.
\ \\

{\bf Aknowledgment : } I would like to thank warmly Karin Melnick for very interesting conversations on the topic. This work was supported by the ANR {\sc Geodycos}.

\end{document}